\documentclass[a4paper, 12pt]{article}

\usepackage[english]{babel}
\usepackage[utf8]{inputenc}
\usepackage[pdftex]{graphicx}
\usepackage[left=1in, right = 1in]{geometry}

\usepackage{authblk}
\usepackage{amsthm,amsmath,amssymb}
\usepackage{algorithm}
\usepackage{algorithmic}
\usepackage{hyperref}
\usepackage{graphicx}
\usepackage{comment}
\usepackage{array}

\newcommand{\R}{\mathbb{R}}

\newcommand{\x}{\mathbf{x}}
\newcommand{\xstar}{\mathbf{x}^\star}
\newcommand{\xhat}{\hat{\mathbf{x}}}
\newcommand{\y}{\mathbf{y}}
\newcommand{\w}{\mathbf{w}}
\newcommand{\argmax}{\operatorname{argmax}}

\newcommand{\abs}[1]{\left| #1 \right|}
\newcommand{\norm}[1]{\left| #1 \right|}
\renewcommand{\brace}[1]{\left\{ #1 \right\}}
\renewcommand{\epsilon}{\varepsilon}

\newcolumntype{L}[1]{>{\raggedright\let\newline\\\arraybackslash\hspace{0pt}}m{#1}}
\newcolumntype{C}[1]{>{\centering\let\newline\\\arraybackslash\hspace{0pt}}m{#1}}
\newcolumntype{R}[1]{>{\raggedleft\let\newline\\\arraybackslash\hspace{0pt}}m{#1}}

\newtheorem{axiom}{Axiom}[section]

\newtheorem{definition}[axiom]{Definition}
\newtheorem{example}[axiom]{Example}

\usepackage{natbib}
 \bibpunct[, ]{(}{)}{,}{a}{}{,}%

\usepackage{tikz}
\usepackage{soul}

\begin{document}
%%%%%%%%%%%%%%%%

\title{Optimising Electric Vehicle Charging Station Placement using Advanced Discrete Choice Models}

\author[1]{Steven Lamontagne}
\author[1]{Margarida Carvalho}
\author[1]{Emma Frejinger}
\author[1]{Bernard Gendron}
\author[2]{Miguel F. Anjos}
\author[3]{Ribal Atallah}

\affil[1]{CIRRELT and D\'epartement d'informatique et de recherche op\'erationnelle, Universit\'e de Montr\'eal} 
\affil[2]{School of Mathematics, University of Edinburgh}
\affil[3]{Institut de Recherche d'Hydro-Qu\'ebec}

%\begin{abstract}We present a new model for finding the optimal placement of electric vehicle charging stations across a multi-period time frame so as to maximise electric vehicle adoption.  Via the use of stochastic discrete choice models and user classes, this work allows for a granular modelling of user attributes and their preferences in regard to charging station characteristics. We adopt a simulation approach and pre-compute error terms for each option available to users for a given number of scenarios. This results in a bilevel optimisation model that is, however, intractable for all but the simplest instances.  {Our major contribution is a reformulation into a maximum covering model, which uses  the pre-computed error terms to calculate the users covered by each charging station. This allows solutions to be found more efficiently than for the bilevel formulation.} The maximum covering formulation remains intractable in some instances, so we propose rolling horizon, greedy, and GRASP heuristics to obtain good quality solutions more efficiently. Extensive computational results are provided, which compare the maximum covering formulation with the current state-of-the-art, both for exact solutions and the heuristic methods.
%{

  %\small	
  %\textbf{\textit{Keywords---}} Electric vehicle charging stations, facility location, integer programming, discrete choice models, maximum covering
%}
%\end{abstract}

\maketitle

\section{Introduction}

In order to meet CO$_2$ emissions goals, major changes at the global level are necessary. This includes the transportation sector, which in 2019 was responsible for 27\% of global CO$_2$ emissions~\citep{IEA2021}  and 30\% of emissions in Canada~\citep{ECCA2021}. One initiative that has been proposed to reduce these emissions is the adoption of electric vehicles (EVs), rather than internal combustion engine vehicles (ICEVs). Depending on the electricity generation mix, this may be effective in reducing emissions~\citep{CPEV2015, Woo2017}. 

Due to the limits of the internal battery, EVs require more frequent refuelling than CVs. While early adopters of EVs may have access to home charging~\citep{Bailey2015}, this can be supplemented with public charging infrastructure to increase accessibility. The availability of the latter has been found to increase EV adoption~\citep[e.g.][]{Coffman2017}. However, this scenario has been noted to suffer from a ``chicken-and-egg'' dilemma~\citep{Anjos2020}, where users require charging infrastructure in order to purchase an EV, while businesses and infrastructure operators have little incentive to install charging stations when there are few users.  To this end, governmental organisations can alleviate this problem by investing in public charging infrastructure, thus allowing users to recharge EVs and encourage EV adoption. 

The problem we examine is that of the decision makers responsible for public charging infrastructure, and the optimal placement of that infrastructure within a city. They have a set of candidate locations for charging stations, and must decide which charging stations to open as well as how many charging outlets to place at each open station. These decisions take place over a long-term planning horizon, and investment in each period is limited by a budget. The decision makers take into account the users who are purchasing a vehicle in each period. Depending on the placement of charging infrastructure, these users may purchase an EV. We assume that the users anticipate the need to recharge an EV, and will only purchase one if they have access to charging infrastructure (either public or at home). Additionally, we assume that more convenient charging (e.g. by a charging station being closer to home, or having more charging outlets) increases the probability that a user will purchase an EV. The goal of the decision makers is to plan the placement of charging infrastructure so that it maximises EV adoption. 

The existing literature for this problem is quite narrow, and has limitations that we address in this work. In the intercity context, the optimisation model of~\cite{Zhang2017} accounts for EV adoption growth depending on the coverage of paths, but their method is not applicable to the intracity case. To the best of our knowledge, in the intracity context, only \cite{Anjos2020} formulate an optimisation model for locating charging stations accounting for EV adoption growth based on the location of charging stations. There, the total number of EVs increases according to a piecewise-linear growth function, which accounts for natural growth. This growth function is applied indiscriminately to all locations.  Additionally, the fraction of users who have access to home charging is assumed to be constant across locations in the network. In reality, the EV growth rate varies by location~\citep{AVEQ2021}, and access to home charging depends on the type of residences~\citep{ICCT2019}.

This work presents modelling, algorithmic, and computational contributions.  For the modelling side, we consider user classes, allowing for parameter values to be considered more specifically for groups of users rather than the entire population. Additionally, within each user class, the optimisation model supports the use of advanced discrete choice models.  {Acknowledging the uncertainty associated with demand distribution, these models assign a probability distribution over available alternatives.} This allows to model heterogeneous user preferences and complex substitution patterns.
%the level of heterogeneity to vary based on the user class. 
The combination of the user classes and the ability to use advanced discrete choice models results in a highly flexible model, which the decision makers can design to suit their specific problem and available data sources. The user behaviour is incorporated in the decision-maker's problem as in \cite{Paneque2020}, leading to a bilevel program. While the bilevel optimisation model presented here can only solve small instances, the maximum covering reformulation allows for instances to be solved significantly faster.  {To the best of our knowledge, such a reformulation has not been presented in the context of the simulation-based approach of \cite{Paneque2020}. We note that this reformulation is effective for our optimisation model, but is also applicable  to more general, uncapacitated bilevel models using the same framework.} On the algorithmic side, we propose three heuristics to solve the optimisation models. 
%the outlines and code for three heuristic methods are applied to our optimisation model. 
First, a rolling horizon method, which was also used in \cite{Anjos2020}. Second, a greedy method~\citep{Church1974}, and third, a GRASP method~\citep{Resende1998}. We extend the latter two to a multi-period setting with sizing decisions.
%from the classical maximum covering problem to work in the multi-period, multi-size maximum covering problem. 
On the computational side, we present extensive experimental comparisons examining: the capabilities of our optimisation model versus the model presented in \cite{Anjos2020},  {the bilevel formulation versus maximum covering formulation,} and the heuristic methods. Where possible, parameter values and user class characteristics are based on real data.

We structure the remainder of the paper as follows: Section~\ref{SectionLiteratureReview} reviews the relevant literature for the problem. Section~\ref{SectionProblemFormulation} describes the framework of the problem,  {as well as briefly presenting the bilevel model and its single-level reformulation. }  Section~\ref{SectionMaximumCovering} gives the maximum covering reformulation for the bilevel model. Section~\ref{SectionHeuristics} discusses various heuristic methods for solving large instance sizes  {of practical interest}. Section~\ref{SectionResults} provides computational results.

\section{Literature Review}
\label{SectionLiteratureReview}
We review three topics of research that relate closely to our work. First, EV charging station location models, which determines the optimal location for charging stations. Second, vehicle acquisition, more specifically relating to EVs, which discusses factors that affect whether users purchase EVs. Third, Maximum Capture Problem with Random Utilities (MCPRU), which examines the optimal placement of facilities in a competitive environment.

On the topic of EV charging station location models, the model proposed by \cite{Anjos2020}, which we refer to as the \emph{Growth Function} (GF) model, is the closest to our work. Whereas most charging station location models are examined with the objective to maximise profits or minimise costs, the GF model is designed solely to encourage electric vehicle adoption. Its use of a combination of node-based and path-based approaches allows for the inclusion of a general highway charging station network while allowing for increased precision within the city. EV adoption and charging stations are linked via capacity constraints, with users willing to recharge at any charging station within a fixed distance of their home.

There is an existing literature for EV charging station placement~\citep[an excellent review is presented in][]{Kadri2020}, however two notable differences separate it from our work. First, in all cases, the users under consideration are those who already own EVs, and are deciding on a charging station to recharge. Second, the objective of the decision maker are different, such as profit maximisation~\citep{Luo2015}, maximising the EV flow that can travel each path given a limited EV range~\citep{Lim2010}, minimising the users' costs for recharging the EV and travel time~\citep{Gonzalez2022}, or maximising the EV charging demand that can be covered~\citep{Frade2011}. We provide a summary of the key characteristics of related works in Table~\ref{TableChargingStationModels}. 

Unlike our work, EV recharging demand at each station or along each path is generally taken to be deterministic. However, in \cite{Luo2015} and \cite{Cui2018} the demand for each station is estimated using the analytic choice probabilities from a nested logit model. The demand model they use considers station characteristics such as the distance to the EV owner and the proximity of the charging station to amenities (such as restaurants and shopping centres). In \cite{Kadri2020}, the evolution each year of the demand for each path is modelled with a  discrete scenario tree. In the tree, the demand for the current year is assumed to be known, as are the transition probabilities for each of the possible states in each subsequent year. This allows for the estimation of the expected demand of each path considering all subsequent years. 

The inclusion of capacity constraints in EV charging station location models is not consistent. Several models (including our work, as well as others in Table~\ref{TableChargingStationModels}) do not consider capacity. For those that do include capacity constraints, the modelling assumptions differ. In \cite{Frade2011} and \cite{Zhang2017}, the demand that may be satisfied by each charging station is limited by the number of outlets installed at that location. In \cite{Luo2015}, the decision maker must ensure that charging stations meet minimum quality of service requirements, including waiting time and service coverage. In \cite{Cui2018}, the decision maker must ensure that the amount of electricity supplied at each charging station is sufficient for the expected amount of demand. In \cite{Gonzalez2022}, the number of recharging sessions that will be required at each station is estimated, and used as a bound for the capacity of the station.

\begin{table}
\small
\begin{tabular}{l R{3cm} R{1.75cm} R{1.5cm} R{1.5cm} R{1.5cm} R{1.5cm}}
\hline
     &  Objective & Model Type & Demand & Time Periods & Capacity & Intracity or Intercity\\
\hline
\cite{Lim2010} & Maximise flow refueled & FRLM & D & Single & No & Inter
\\
\hline
\cite{Frade2011} & Maximise coverage of EV charging demand & Maximum Covering & D & Single & Yes & Intra
\\
\hline
\cite{Shukla2011} & Maximise flow intercepted & FILM & D & Single & No & Intra \\
\hline
\cite{Capar2013} & Maximise flow refueled & FRLM & D & Single & No & Inter
\\
\hline
\cite{Luo2015} & Maximise profit of decision maker & MINLP & S & Multi & Yes & Intra
\\
\hline
\cite{Zhang2017} & Maximise flow refueled & FRLM & D & Multi & Yes & Inter
\\
\hline
\cite{Cui2018} & Minimise cost of decision maker & MINLP & S & Multi & Yes & Intra
\\
\hline
\cite{BadriKoohi2019} & Minimise costs, p-median, flow interception (weighted sum) & MILP & D & Single & No & Intra
\\
\hline
\cite{Kadri2020} & Maximise flow refueled & FRLM & S & Multi & No & Inter
\\
\hline
\cite{Gonzalez2022} & Minimise users' travel cost & Bilevel & D & Single & Yes & Intra\\
\hline
\multicolumn{7}{l}{Legend: FRLM = Flow Refueling Location Model, FILM = Flow Interception Location Model,} \\
\multicolumn{7}{c}{MINLP = Mixed Integer Non-Linear Problem, MILP = Mixed Integer Linear Problem,} \\
\multicolumn{7}{c}{D = Deterministic, S = Stochastic.}
\end{tabular}
\caption{Summary of key characteristics of EV charging station location models.}
\label{TableChargingStationModels}
\end{table}

%In terms of models which specifically consider EV growth based on charging station placement, \cite{Anjos2020} and \cite{Zhang2017} both utilise capacity constraints. In the former, the capacity forms an upper bound for the growth of EVs, forbidding users from purchasing EVs once the capacity has been exceeded. In the latter, the capacity is included as one of the factors in the EV growth, increasing the number of EVs beyond the natural growth based on the network coverage.  

The vehicle acquisition problem examines the attributes that increase the likelihood of purchasing EVs, and are useful for determining exogenous attributes for the demand model and defining user classes. Importantly, the availability of electric vehicle charging stations increases the likelihood of users purchasing an EV; This availability (also referred to as `fuel availability'), was a significant factor in the decision to purchase an EV in most articles examining the subject~\citep{
Achtnicht2012, Javid2016, Hackbarth2013, Ziegler2012,
Rezvani2015, Coffman2017}. This finding is not unanimous, however, as \cite{Bailey2015} propose that a predisposition to EVs makes the users more likely to notice existing charging infrastructure, not the converse. The conclusions of \cite{Bailey2015} are also present in \cite{Axsen2015}, an article based on the same dataset and project. Several studies focus on assessing preferences in regard to various types of attributes~\citep{Hidrue2011, Achtnicht2012, Ziegler2012, Hackbarth2013, Axsen2015, Bailey2015, Javid2016}. For reviews of EV acquisition models, we refer to \cite{Rezvani2015}, \cite{Javid2016}, and \cite{Coffman2017}. 

In this work, we formulate the problem as a MCPRU, generally attributed to \cite{Benati1999} and \cite{Benati2002}. In the MCPRU, a company is looking to place facilities in an environment where competitors have existing facilities. A \emph{Random Utility Maximisation} (RUM) discrete choice model is used to  {assign choice probabilities of users}, based on the set of facilities available to them. The company aims to place the new facilities to maximise the market share that the new facilities capture. In general, the MCPRU focuses on two key attributes of facilities in order to determine which one users choose to patronise: distance to the users, and \emph{facility attractiveness}~\citep{Berman2014}. It is possible to consider market expansion in the MCPRU, where the placement of facilities attracts customers that were not originally in the market~\citep[see, e.g.,][]{Aboolian2007}. Since a vehicle (EV or otherwise) is an expensive purchase, it is assumed that the additions to the charging infrastructure are not sufficient to attract users who were not planning on purchasing a vehicle at all in the given period. We thus do not include market expansion in our optimisation model.

Almost all work on the MCPRU has been done in the context of using the multinomial logit (MNL) model for characterising user behaviours, due to the existence of an analytic formula for calculating choice probabilities. For recent examples, the linear formulation in \cite{Freire2016} \citep[which improves on the linear formulation in][]{Haase2009} and the outer approximation and submodular cut methods in \cite{Ljubic2018} both rely on the logit choice probabilities. Exceptions to the use of MNL model include the work of \cite{Mai2020}, who express a mixed logit choice probability as a sum of probabilities of MNL models (though they do not conduct any computational experiments using this method), and the work of \cite{Dam2021}, who use a method that makes use of the submodularity of the objective function for all discrete choice models in the Generalised Extreme Value (GEV) family (of which MNL and mixed logit are included). Since we do not assume the use of MNL or the GEV family (in addition to the added complexity of multi-period and sizing considerations), the methods proposed in these works cannot be applied directly to our model.

While it can be used more generally, the simulation-based approach in \cite{Paneque2020} can be applied to the MCPRU. Rather than embedding an analytic probability that users will select an alternative into the objective function, this approach instead generates error terms for each alternative for a given number of \emph{scenarios}. These error terms, pre-computed  and given to the optimisation model as an input, allow for the utility to be calculated for each alternative and each scenario. Using the RUM principle, users then select, in each scenario, the alternative which has the highest utility. Alternatives which are not available to users (such as, in the case of the MCPRU,  facilities which are closed) are set to a lower bound, ensuring they are not selected. This approach supports the use of any discrete choice model rather than being limited to the MNL model (or even the GEV family). However, it is computationally difficult to solve in larger instances.

 {
We emphasise that the simulation-based approach of \cite{Paneque2020} differs from the linearisation typically presented for the MCPRU \citep[e.g. ][]{Benati2002, Haase2009, Haase2014, Freire2016}. In the latter case, the analytic choice probabilities for each alternative are linearised by setting the variable coefficients appropriately. However, this linearisation can only be done when analytic choice probabilities are known for the discrete choice model (e.g. MNL). On the contrary, the simulation-based approach assumes that the utility functions are linear in the decision variables and works directly in the space of utilities. Under the RUM principle, this allows for linearisation of the sample average approximation of the choice probabilities (which can be derived for any discrete choice model for which error terms can be drawn). As a consequence, the choice probabilities are approximations of their exact form, but also renders the approach suitable for any discrete choice model. This includes, notably, discrete choice models for which the choice probabilities are typically estimated via sample average approximation \citep[e.g. mixed logit][]{Train2002}.
}

Our work aims to bridge the gap between the EV charging station placement models, EV acquisition models, and the MCPRU. Most EV charging station placement models are designed to benefit existing EV owners, instead of maximising EV adoption. The existing works which include EV growth based on the placement of EV charging stations, namely the models of \cite{Zhang2017} and \cite{Anjos2020}, use an aggregate approach for EV growth. However, the literature on EV acquisition shows that significant factors for determining EV acquisition are individual characteristics of the users (e.g. income, education level, access to home charging). This highlights the benefit of user classes in our optimisation model, as they allow for these individual characteristics to be considered. Moreover, we note the heterogeneity within each preference class (which function equivalently as a user class in our work) in the latent class models of \cite{Hidrue2011} and \cite{Axsen2015}. This suggests the use of a discrete choice model for demand modelling, and naturally leads to a MCPRU formulation for the optimisation model. As is the case for the MCPRU, in addition to user characteristics,  characteristics from the charging station also impact the choice of users. While the specification of the discrete choice model for the users is outside the scope of this work, we note  {in Appendix~\ref{AppendixBilevelModel}} characteristics of the demand modelling which makes the MNL model less suitable. Because of this, we use the simulation-based approach of \cite{Paneque2020}.

\section{Problem Formulation and Modelling}
\label{SectionProblemFormulation}

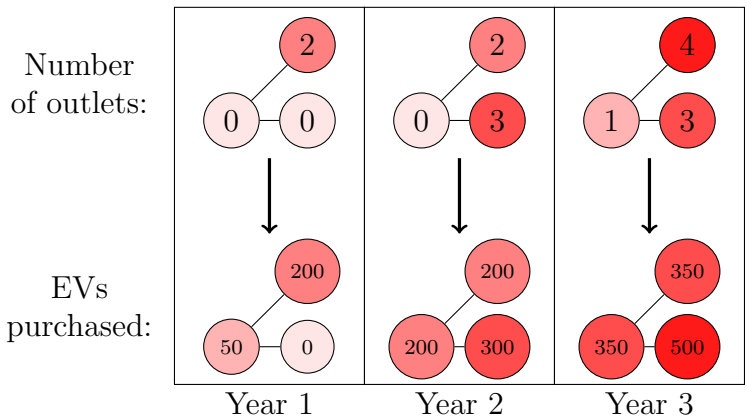
\begin{figure}
\centering
% \documentclass{article}

% \usepackage{tikz}
% \usetikzlibrary{external}
% \tikzexternalize[prefix=Images/]
% \begin{document}

\begin{tikzpicture}[ZeroOutlet/.style={circle, draw = black, fill=red!10, minimum size=7mm, align = center},
OneOutlet/.style={circle, draw = black, fill=red!30, minimum size=7mm, align = center},
TwoOutlet/.style={circle, draw = black, fill=red!50, minimum size=7mm, align = center},
ThreeOutlet/.style={circle, draw = black, fill=red!70, minimum size=7mm, align = center},
FourOutlet/.style={circle, draw = black, fill=red!90, minimum size=7mm, align = center},
Text/.style={text width = 2cm, align = center},
]

\node[Text] (Stations0) at (-2, 0.5) {Number of outlets:
};
\node[Text] (EVs0) at (-2, -2.5) {EVs purchased:
};

\node[ZeroOutlet] (s00) at (0,0) {
0};
\node[ZeroOutlet] (s01) at (1,0) {
0};
\node[TwoOutlet] (s02) at (1,1) {
2};
\draw[-] (s00) -- (s01);
\draw[-] (s00) -- (s02);

\node[OneOutlet] (e00) at (0,-3) {
\scriptsize{50}};
\node[ZeroOutlet] (e01) at (1,-3) {
\scriptsize{0}};
\node[TwoOutlet] (e02) at (1,-2) {
\scriptsize{200}};
\draw[-] (e00) -- (e01);
\draw[-] (e00) -- (e02);

\draw[black] (-0.75,1.5) rectangle (1.75,-3.5);
\node (Year0) at (0.5,-3.75) {Year 1};
\draw[->, very thick] (0.5,-0.5) -- (0.5, -1.5);

\node[ZeroOutlet] (s10) at (2.5,0) {
0};
\node[ThreeOutlet] (s11) at (3.5,0) {
3};
\node[TwoOutlet] (s12) at (3.5,1) {
2};
\draw[-] (s10) -- (s11);
\draw[-] (s10) -- (s12);

\node[TwoOutlet] (e10) at (2.5,-3) {
\scriptsize{200}};
\node[ThreeOutlet] (e11) at (3.5,-3) {
\scriptsize{300}};
\node[TwoOutlet] (e12) at (3.5,-2) {
\scriptsize{200}};
\draw[-] (e10) -- (e11);
\draw[-] (e10) -- (e12);

\draw[black] (1.75,1.5) rectangle (4.25,-3.5);
\node (Year1) at (3,-3.75) {Year 2};
\draw[->, very thick] (3,-0.5) -- (3, -1.5);

\node[OneOutlet] (s20) at (5,0) {
1};
\node[ThreeOutlet] (s21) at (6,0) {
3};
\node[FourOutlet] (s22) at (6,1) {
4};
\draw[-] (s20) -- (s21);
\draw[-] (s20) -- (s22);

\node[ThreeOutlet] (e20) at (5,-3) {
\scriptsize{350}};
\node[FourOutlet] (e21) at (6,-3) {
\scriptsize{500}};
\node[ThreeOutlet] (e22) at (6,-2) {
\scriptsize{350}};
\draw[-] (e20) -- (e21);
\draw[-] (e20) -- (e22);

\draw[black] (4.25,1.5) rectangle (6.75,-3.5);
\node (Year1) at (5.5,-3.75) {Year 3};
\draw[->, very thick] (5.5,-0.5) -- (5.5, -1.5);
\end{tikzpicture}

% \end{document}
\caption{The decision maker decides where to place charging infrastructure (stations and number of outlets) on the simple example network. In response, a subset of users decide to purchase EVs.}
\label{FigureProblemFormulation}
\end{figure}

In our problem, we consider the placement of charging infrastructure in a city, wherein there are two parties with different motives: a decision maker placing the infrastructure, and the users choosing whether to purchase an EV (depending on the charging infrastructure). An illustration of this process is given in Figure~\ref{FigureProblemFormulation}.  

\subsection{Decision Maker}
The decision maker is planning charging infrastructure over a planning horizon with $T$ periods (e.g., seasons or years) indexed by $1 \leq t \leq T$.

The decision maker has a candidate set of locations $M$ where charging stations $j \in M$ may be installed or expanded with additional outlets to a maximum number $m_j$. In each period $t$, there is a budget $B^t$, which limits the total investment (both opening stations and installing outlets).  {Let $c_{jk}^t$ denote the cost to go from $k-1$ to $k$ outlets at charging station $j \in M$ in period $t$. Note that this cost concept is flexible, and accounts for additional costs to be incorporated at certain thresholds if significant infrastructure upgrades are required. Notably, this includes initial infrastructure installation for $k=1$ but can also include, e.g., increasing electric power capacity at the station.}

 {Let $x_{jk}^{t}$ be a binary variable indicating if station $j \in M$ has \emph{at least} $k$ outlets in period $t$ (with $1 \leq k \leq m_j$).}  We denote $\x = \brace{x_{jk}^t}, j \in M, 1 \leq k \leq {m_j}, 1 \leq t \leq T$.  {Let $x_{jk}^0$ denote the initial state of each charging station $j$.}

\subsection{Users}
We consider users who are are planning on purchasing a vehicle and, depending on public charging infrastructure, may choose to purchase an EV. These users are composed of a set of user classes $N$. The population size for each user class is given by $N_i^t, i \in N, 1 \leq t \leq T$.  {These parameters provide a flexible modelling. For instance, they can correspond to the number of potential new EVs in a class or be configured to prioritise early adopters by decreasing the population size in later time periods. This would have an effect similar to the classical approach of discounting future time periods.}

Each user considers whether they have a primary option for recharging an EV. If such an option is available and sufficient for their needs, they decide to purchase an EV. This decision is modelled via a discrete choice model, with users in user class $i \in N$ selecting an alternative $j$ in choice set  {$\mathcal{C}_i^t(\x)$}. The inclusion of  {$\x$} in the notation for the choice set emphasises the dependence on the set of open stations, since stations  {without at least one outlet} cannot be used to recharge an EV. The choice set   {$\mathcal{C}_i^t(\x)$} contains a subset of the public charging infrastructure, home charging (if available to  user class $i$), and an opt-out. The latter corresponds to the alternative to not purchase an EV, while all other alternatives indicate the user purchases one. In what follows, we denote the opt-out alternative as $j=0$.
 { We note that the use of the index $j$ for both alternatives and stations is deliberate because, as we discuss in Section~\ref{SectionSimulationBasedModel}, the charging stations are the alternatives for the most part.}

 { At the strategic-level planning,} we assume that each charging station $j \in M$ that is opened is \emph{uncapacitated}, in the sense that there is no limit to the number of users that may select it. However, we assume that the likelihood of users selecting a given station increases with the number of outlets. The rationale being that the higher the number of outlets, the more likely one may be available when required. This allows us to implicitly take into account that users perceive the capacity as finite.

\subsection{Simulation-based Model}
\label{SectionSimulationBasedModel}
 {We adopt the simulation-based method of \cite{Paneque2020} to the users' decision to purchase an EV. More specifically, for each user class $i\in N$, let $R_i$ be the number of scenarios used to approximate the choice probabilities. }
We assume that the  {decision to purchase an EV} is separable for each user class, each period, and each scenario, meaning that there is no interaction among them. Therefore, in what follows, we concentrate on detailing  {a given triplet $(t,i,r) \in P =\{(t',i',r'): 1\leq t' \leq T, i' \in N, 1\leq r' \leq R_{i'}\}$}.  For each alternative $j \in  {\mathcal{C}_i^t (x)}$, let $u_{ji}^{rt}$ be the \emph{simulated utility}.

We previously remarked that the choice set  {$\mathcal{C}_i^t(\x)$} depends on which stations are open. For notational simplicity, we define the two choice sets: $C_i^{0t}$ and $C_i^{1t}$.  The set $C_i^{0t}$ includes all alternatives exogenous to the optimisation model (i.e. unaffected by decision variables), for which the availability is not under the control of the decision maker. This always includes the opt-out alternative, but may also, for example, include the alternative for home charging. The simulated utility $u_{ji}^{rt}$, for $j \in C_i^{0t}$,  is thus given by
\begin{equation}
\label{UtilityCalc_Exo}
u_{ji}^{rt}=\kappa_{ji}^t +\epsilon_{ji}^{rt},
\end{equation}
 where the parameter $\kappa_{ji}^t$ is the alternative-specific constant. In general, the alternative-specific constants can be any function which does not depend on the decision variables.  {Note that we assume only one opt-out alternative as multiple options can be replaced by a single alternative given by the maximum utility across all opt-out alternatives.}

The set $C_i^{1t}$ includes the alternatives for charging stations which have a non-zero probability of being chosen (e.g. sufficiently close to be considered). If station $j \in C_i^{1t}$ is open, we consider the utility $u_{ji}^{rt}$ to be a function of the number of charging outlets available at that station as well as the alternative-specific constant. This leads to the formulation:
%\begin{equation*}
$  
u_{ji}^{rt}= \sum_{k=1}^{m_j} \beta_{jik}^t x_{jk}^t+\kappa_{ji}^t +\epsilon_{ji}^{rt}.$
%\end{equation*}
 {We recall that $m_j$ is the maximum number of outlets at station $j$. The parameter $\beta_{jik}^t$ is the incremental benefit of going from $k-1$ to $k$ outlets.} The model does not impose any restrictions on the form of the utility coefficients, however,  {we assume that the coefficients are non-negative}. As before, the alternative-specific constant can be any function which does not depend on the decision variables.
 {If the station is closed, the utility is set to a lower bound $\underline{a}_i^{t}$ to ensure that it will not be selected. }

\subsection{Single-Level Model}\label{subsec:singlelevel}

 {Given $\x$ and a triplet $(t,i,r)$, a user selects the alternative maximising the simulated utility, i.e., $k \in \arg\max_{j \in \mathcal{C}_i^t(\x)} u^{rt}_{ji}$. If we now integrate the goal of the decision maker, which is to maximise the purchase of EVs, we naturally obtain a bilevel formulation. Following the simulation-based method of \cite{Paneque2020}, we can reduce the bilevel formulation to a single-level form. Indeed, the Single-Level (SL) model of the problem is given by}
\begin{subequations}
    \begin{alignat}{3}
\operatorname{Minimise} & \sum_{(t,i,r) \in P} \frac{N_i^{t}}{R_i} w_{0i}^{rt}, 
\label{ModelSinglelevel:Objective}
\\
\text{subject to }& \sum_{j \in M}\sum_{k=1}^{m_j} c_{jk}^ t\left(x_{jk}^ t-x_{jk}^{t-1} \right)  \leq B^t, && 1 \leq t \leq T. 
\label{ModelSinglelevel:Budget} 
\\
& x_{jk}^t \leq x_{jk-1}^t,  && 1 \leq t \leq T, j \in M, 1 \leq k \leq m_j,
\label{ModelSinglelevel:AtLeastKOutlets} 
\\
& x_{jk}^t \geq x_{jk}^{t-1}, && 1 \leq t \leq T, j \in M, 1 \leq k \leq m_j, 
\label{ModelSinglelevel:KeepOutlets}
\\
& u_{ji}^{rt}=\kappa_{ji}^t +\epsilon_{ji}^{rt},  && (t,i,r) \in P, j \in C_i^{0t}, 
\label{ModelSinglelevel:UtilityCalc_Exo}
\\
& u_{ji}^{rt} \geq  {\underline{a}_i^t},  && (t,i,r) \in P,  j \in C_i^{1t},  
\label{ModelSinglelevel:DiscountClosed1}
\\
& u_{ji}^{rt} \leq  {\underline{a}_i^t} + \nu_{ji}^{rt} x_{j1}^t, &&  (t,i,r) \in P,  j \in C_i^{1t},    
\label{ModelSinglelevel:DiscountClosed2}
\\
& u_{ji}^{rt} \geq \sum_{k=1}^{m_j} \beta_{jik}^t x_{jk}^t+\kappa_{ji}^t +\epsilon_{ji}^{rt} - \nu_{ji} ^{rt} \left( 1 - x_{j1}^t \right), &&  (t,i,r) \in P,  j \in C_i^{1t},   \label{ModelSinglelevel:DiscountOpen1}
\\
& u_{ji}^{rt} \leq \sum_{k=1}^{m_j} \beta_{jik}^t x_{jk}^t+\kappa_{ji}^t +\epsilon_{ji}^{rt}, &&  (t,i,r) \in P,  j \in C_i^{1t},   
\label{ModelSinglelevel:DiscountOpen2}
\\
& u_{ji}^{rt}-\alpha_{i}^{rt}+\left(1-w_{ji}^{rt}\right)\mu_{ji}^{rt}  \geq 0, && (t,i,r) \in P, j \in C_i^{0t} \cup C_i^{1t}, 
\label{ModelSinglelevel:LowerLevelSolutionStart}
\\
& \sum_{j \in C_i^{0t}} w_{ji}^{rt}+\sum_{j \in C_i^{1t}} w_{ji}^{rt}=1, && (t,i,r) \in P,
\\
& \alpha_{i}^{rt} \geq u_{ji}^{rt},  && (t,i,r) \in P,  j \in C_i^{0t} \cup C_i^{1t},
\label{ModelSinglelevel:LowerLevelSolutionEnd}
\\
& w_{ji}^{rt} \in \{0,1\},  && (t,i,r) \in P, j \in C_i^{0t} \cup C_i^{1t}, \notag
\\
& \alpha_{i}^{rt} \in \mathbb{R}, && (t,i,r) \in P,
\notag
\\
& u_{ji}^{rt} \in \mathbb{R}, && (t,i,r) \in P, j \in C_i^{0t} \cup C_i^{1t},
\notag
\\
& x_{jk}^t \in \{0,1\}, && 1 \leq t \leq T, j \in M, 1 \leq k \leq m_j. 
\notag
\end{alignat}
\label{ModelSinglelevel}
\end{subequations}
The objective function~\eqref{ModelSinglelevel:Objective} minimises the number of users who do not purchase an EV (or, equivalently, maximises those who do) with the auxiliary variables $w_{ji}^{rt}$ representing the selection of alternative $j$ by the triplet $(t,i,r)$. Recall that $w_{0i}^{rt}$ corresponds to the opt-out alternative, and indicates that the triplet $(t,i,r)$ does not purchase an EV.  Budget constraints~\eqref{ModelSinglelevel:Budget} are given for each time period $t$. It would be possible to supplement the per-period budget with an overall budget, as was done in \cite{Anjos2020}.  {Constraints~\eqref{ModelSinglelevel:AtLeastKOutlets} impose that to have at least $k$ outlets, we must have at least $k-1$ outlets.} Constraints~\eqref{ModelSinglelevel:KeepOutlets} forbid us from removing outlets once installed. Constraints~\eqref{ModelSinglelevel:UtilityCalc_Exo} set the utility for exogenous alternatives (opt-out, home charging, etc.). Constraints~\eqref{ModelSinglelevel:DiscountClosed1}-\eqref{ModelSinglelevel:DiscountOpen2} set the utility for endogenous alternatives (charging stations). Constraints~\eqref{ModelSinglelevel:LowerLevelSolutionStart}-\eqref{ModelSinglelevel:LowerLevelSolutionEnd} ensure that the alternative with the highest utility is selected for each triplet $(t,i,r)$.  We present the details around the model in Appendix~\ref{AppendixBilevelModel}.

We note that Constraints~\eqref{ModelSinglelevel:UtilityCalc_Exo} set the utility for all exogenous alternatives in $C_i^{0t}$, which includes opt-out and may include home charging.  {However, in practice, it is better to preprocess any triplets which have access to home charging: if the utility associated with home charging is lower than that of the opt-out, then home charging can be ignored (as it will never be selected). If it is higher than for the opt-out, then the associated opt-out choice variable $w_{0i}^{rt}$ can be fixed to 0, as we can guarantee that it will not be selected. As such, the set $C_i^{0t}$ can be reduced to only contain the opt-out alternative. }

\section{Maximum Covering Model}
\label{SectionMaximumCovering}
While the SL model~\eqref{ModelSinglelevel} is a Mixed-Integer Linear Programming (MILP) optimisation problem that can be given directly to a general purpose solver, large-scale instances can be hard to solve. This is due to the Big-M constraints for the  utility~\eqref{ModelSinglelevel:DiscountClosed1}-\eqref{ModelSinglelevel:DiscountOpen2} and the linearised lower-level problem~\eqref{ModelSinglelevel:LowerLevelSolutionStart}-\eqref{ModelSinglelevel:LowerLevelSolutionEnd}, as well as both sets of binary variables $\x$ and $\w$.  {While the Big-M values are tight given the bounds, as we will see in Section~\ref{SectionResults},} the model is intractable for all but the simplest of instances. For this reason, we propose to reformulate the problem into a maximum covering problem using the pre-computed error terms $\epsilon_{ji}^{rt}$.

For $j \in C_i^{1t}$, define
%\begin{equation*}
$u_{jik}^{rt} = {\sum_{k'=1}^k \beta_{jik'}^t +\kappa_{ji}^t +\epsilon_{ji}^t}$.
%\end{equation*}
In other words, $u_{jik}^{rt}$ is the utility for the triplet $(t, i,r)$ for charging station $j$ having $k$ outlets in period $t$. This is equivalent to the utility $u_{ji}^{rt}$ with  {$x_{jk'}^t=1 \text{ if } k' \leq k \text{ and } x_{jk'}^t =0 \text{ if } k' > k$.}
\begin{definition}
A charging station $j$ with $k$ charging outlets \emph{covers} the triplet $(t, i, r)$ if the following conditions hold: \emph{(i)} $k \geq 1$, \emph{(ii)} $j \in C_i^{1t}$, and \emph{(iii)} $u_{jik}^{rt} \geq  u_{0i}^{rt}$, where $u_{0i}^{rt}$ represents the opt-out utility for triplet $(t, i , r)$. We say that $(t,i,r)$ is covered by $\x$ if $\exists j \in M, \exists k \in \{1, \dots , m_j\}$ such that $x_{jk}^t = 1$ and charging station $j$ with $k$ charging outlets covers $(t,i,r)$.
\end{definition}

Intuitively, a charging station covers a triplet if the station has at least one outlet (and is thus open), it is available to and considered by the user class in question, and the charging station is a better option than the opt-out for that triplet. %We give an example of the calculation of a covering in Example~\ref{ExampleUtilitiesCovering}.  

\begin{example}
\label{ExampleUtilitiesCovering}
We consider a given scenario $r$, user class $i$, and period $t$. In Figure~\ref{example:cover}, we see the utilities $u_{ji}^{rt}$ for each option $j \in C_i^{0t} \cup C_i^{1t}$. For each option, the error terms and alternative-specific constants are all pre-computed, which defines the values for utility with no charging outlets. The opt-out utility does not depend on the number of outlets, and so it is constant. Stations 1 and 3 only cover the triplet $(t,i,r)$ if they have at least one and four outlets, respectively. Station 2 is not able to cover $(t,i,r)$ even with six outlets.

\begin{figure}[!ht]
\centering{
\includegraphics[width=0.6\textwidth]{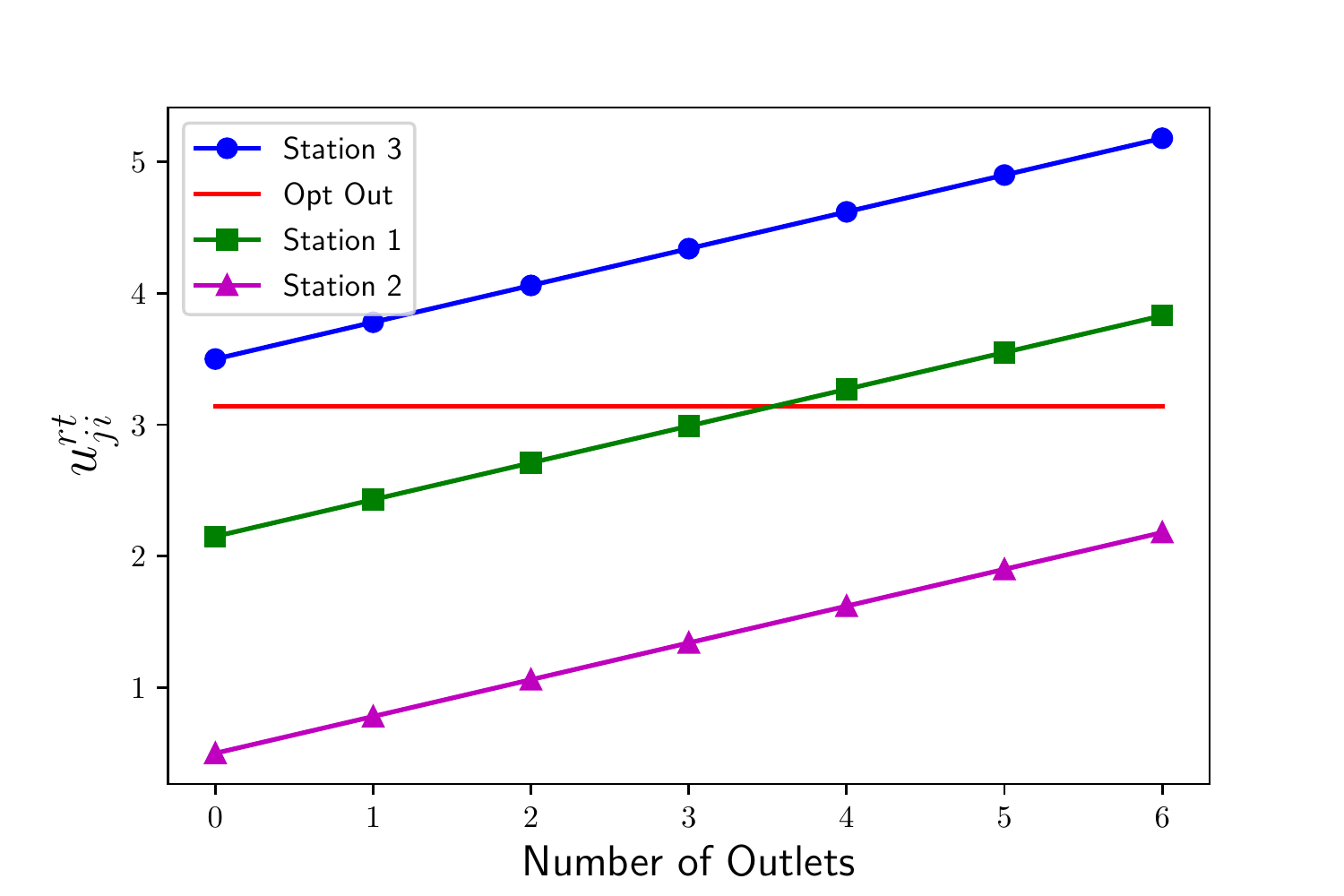}
\caption{Utilities for user class $i$, under scenario $r$, at period $t$.}\label{example:cover}
}
\end{figure}

\end{example}

Next, we define $a_{jik}^{rt} = 1$ if station $j$ with $k$ outlets covers $(t,i,r)$, and $0$ otherwise. This vector $a$ of parameters can be pre-computed after the error terms have been calculated. In the maximum covering formulation, we keep the  {variables $\x$}. We then define the covering decision variables $w_i^{rt} = 1$ if $(t,i,r)$ is covered by $\x$, and $0$ otherwise.

The Maximum Covering (MC) model is then
 {
\begin{subequations}
\begin{alignat}{3}
\operatorname{Maximise} & \sum_{(t,i,r) \in P} \frac{N_i^{t}}{R_i} w_{i}^{rt}, 
\label{ModelMaximumCover:Objective} 
\\
\text{subject to }
&  \eqref{ModelSinglelevel:Budget} -\eqref{ModelSinglelevel:KeepOutlets} 
\notag
\\
& \sum_{j \in M} \sum_{k=1}^{m_j} a_{jik}^{rt} x_{jk}^t \geq w_i^{rt}, && (t,i,r) \in P, 
\label{ModelMaximumCover:Covering}
\\
& x_{jk}^t \in \brace{0, 1}, && 1 \leq t \leq T, j \in M, 1 \leq k \leq m_j,  \notag 
\\ 
& w_{i}^{rt} \in \left[0,1\right], && (t,i,r) \in P. \notag 
\end{alignat}
\label{ModelMaximumCover}
\end{subequations}
}
The objective function~\eqref{ModelMaximumCover:Objective}   is the maximisation equivalent to the objective in the  {SL} model~\eqref{ModelSinglelevel:Objective}. Constraints~\eqref{ModelMaximumCover:Covering} model whether the triplet $(t, i, r)$ is covered by a given $\x$. We note that due to the direction of the optimisation, these inequalities are satisfied with equality at an optimum. Since these constraints are the only ones on the variables $w_i^{rt}$, it is well-known that their integrality can be relaxed~ \citep{Murray2016}. These constraints replace the utility constraints for choice set $C_i^{0t}$~\eqref{ModelSinglelevel:UtilityCalc_Exo}, the utility constraints for choice set $C_i^{1t}$~\eqref{ModelSinglelevel:DiscountClosed1}-\eqref{ModelSinglelevel:DiscountOpen2} and the  {maximisation of the users' utilities given by} (lower-level problem) constraints~\eqref{ModelSinglelevel:LowerLevelSolutionStart}-\eqref{ModelSinglelevel:LowerLevelSolutionEnd}. They thus eliminate all Big-M constraints from the model.  {It is important to note that the maximum covering reformulation is possible because we consider the case where stations are uncapacitated. However, if we wished to include capacity constraints, this would be straightforward for the SL model~\eqref{ModelSinglelevel} \citep[as was done in ][]{Paneque2020}. 
} 
\section{Heuristic Methods}
\label{SectionHeuristics}

While the reformulation from the SL model~\eqref{ModelSinglelevel} to the MC model~\eqref{ModelMaximumCover} significantly improves tractability, the latter is still unable to solve larger instances, as shown in Section~\ref{SectionResults}. We propose three heuristic methods for solving the MC model~\eqref{ModelMaximumCover}, including rolling horizon, greedy, and Greedy Randomised Adaptive Search Procedure (GRASP).

We note that, for the greedy and GRASP methods, it is required to calculate the quality of a candidate solution $\x$. We denote it by
 {
\begin{equation}
f\left( \x \right) = \sum_{(t,i,r) \in P} \frac{N_i^{t}}{R_i} \min \left(1, \sum_{j \in M} \sum_{k=1}^{m_j} a_{jik}^{rt} x_{jk}^t\right). \label{HeuristicSolutionQuality}
\end{equation}
}

\subsection{Rolling Horizon Heuristic}
 {Several variants of the rolling horizon heuristic exist, which form a standard approach in a multi-period setting.
The basic version treats each period independently.} We use this approach as a baseline. 
More precisely, we solve the MC model~\eqref{ModelMaximumCover} 
  for one period at a time, $1 \leq t \leq T$.  Given the potential difficulty of solving the MC model~\eqref{ModelMaximumCover} even when restricted to one period, a time limit is added. The best incumbent solution found within this time limit is returned by the heuristic.
  
   {An alternate version of the rolling horizon heuristic was tested, where for each time period $1 \leq t \leq T$  the variables $x_{jk}^t$ are discrete, but taken as continuous for $x_{jk}^{t'}, t' >t$. The variables $x_{jk}^t$ are then fixed, and we repeat the procedure for $t = t+1$. However, this method did not perform well in our application, due to the increased time required to solve each time period. Thus, we do not discuss this heuristic further in the rest of the paper.}

\subsection{Greedy Heuristics}
\label{SectionGreedyOne}
 The maximum covering problem is generally attributed to \cite{Church1974}, where they also present a simple greedy heuristic. The algorithm presented here is a natural extension, where we iteratively place outlets one at a time, selecting the one which  {increases the number of new EVs the most. }
 
  {
 More specifically, when considering the number of new EVs covered by a new outlet,} there are two possible search modes. In the \emph{myopic} search mode, only new EVs from the current period are counted. In the \emph{hyperoptic} search mode,  {since any outlets placed in the current time period cannot be removed, the EVs in future time periods that result from these outlets are also counted. Namely, we consider the following \emph{score} functions when examining a candidate solution $\tilde{x}$ in time period $t'$:
\begin{align*}
f_m\left( \tilde{\x}, t' \right) & = \sum_{i \in N}\sum_{r=1}^{R_i}\frac{N_i^{t'}}{R_i} \min \left(1, \sum_{j \in M} \sum_{k=1}^{m_j} a_{jik}^{rt'} \tilde{x}_{jk}^{t'}\right),
\\
f_h\left( \tilde{\x}, t' \right) &= \sum_{t =t'}^T \sum_{i \in N}\sum_{r=1}^{R_i}\frac{N_i^{t}}{R_i} \min \left(1, \sum_{j \in M} \sum_{k=1}^{m_j} a_{jik}^{rt} \tilde{x}_{jk}^{t}\right),
\end{align*}
with $f_m$ and $f_h$ being (respectively) for the myopic and hyperoptic search modes. 

Independently of the search mode, each greedy heuristic starts with the trivial, zero solution (which is always feasible). It then iterates over each time period. Given a time period, for each station, the heuristic determines whether adding an outlet is feasible. If so, it computes the number of new EVs covered accordingly with the score function for the search mode ($f_m$ or $f_h$). The outlet which maximises the value of the score function is then selected, and the candidate solution updated. The iteration over stations for the given period continues until no outlets can feasibly be added (either due to the budget or the maximum number of outlets at each station), or the total number of EVs does not increase after adding a station. The latter is possible if all triplets covered by a possible outlet are already covered in the incumbent solution. The heuristic then proceeds to the next time period and repeats the process. 

} %End of  

\begin{comment}
\begin{algorithm}
{
\begin{algorithmic}[1]
\STATE {Initialise: $\xstar \leftarrow \x^0$, $t \leftarrow 1 $}
\WHILE {$t < T$}
\STATE {Find the station $j$ at which to place an outlet, which remains feasible and increases the objective value the most \label{AlgoGOFindStation}}
\STATE {Set $\x$ as $\xstar$ with one additional outlet at station $j$ \label{AlgoGOUpdateSolutionStart}}
\IF {$f\left( \x \right) > f\left( \xstar \right)$}
\STATE { $\xstar \leftarrow \x$ \label{AlgoGOUpdateSolutionEnd}}
\STATE {Go to line \ref{AlgoGOFindStation}}
\ELSE {\label{AlgoGOIncrementYearStart}}
\STATE {$t \leftarrow t+1$}
\ENDIF{\label{AlgoGOIncrementYearEnd}}
\ENDWHILE
\RETURN $\xstar$
\end{algorithmic} 
}
\caption{Greedy Heuristic} \label{AlgorithmGreedy}
\end{algorithm} 
\end{comment}

\subsection{GRASP}
GRASP is characterised by two phases. The first utilises a greedy approach to generate solutions, where each greedy move is randomly selected from all possible moves that result in a solution with an objective value that is within a factor $\alpha$ of the optimal greedy move. The second phase performs a local search procedure on each solution from the first phase.  This method was first applied to the maximum covering problem by \cite{Resende1998}. However, the algorithm described in this work  makes several adaptations in order to support the multi-period and sizing considerations.

\paragraph{First phase}
The solution construction phase is similar to the Greedy algorithm, with the addition of the parameter $\alpha \in [0,1]$. We randomly select one outlet to place from amongst the set of outlets that result in an increase within $100 \times \alpha$ \% of the best possible outlet. The myopic and hyperoptic search modes also apply in the GRASP method. 

\paragraph{Solution Filtering}
In our method, the second phase of GRASP takes considerably longer than the first phase. As such, it is beneficial to filter out unpromising solutions early on. In a method proposed in \cite{Resende2018}, we examine whether the local search method, when applied to a candidate solution, is likely to result in a better objective value than our incumbent solution. For a given number of candidate solutions, we examine the relative increase to the objective function before and after applying the local search. For all subsequent candidate solutions, we then estimate the maximum objective value from the second phase by multiplying the objective value of the candidate solution by the maximum relative increase that was observed. If this results in an objective value that is less than our current best, we filter the candidate solution and do not start the second phase.

\paragraph{Second Phase}
When examining a candidate solution $\xhat$, we consider the three following moves:
\begin{description}
\item[Add:] If the budget permits, we add an outlet to the given station $j$ in period $t$. To ensure feasibility of the solution, we also increase the number of outlets for all subsequent periods $t+1, \dotsc, T$ for station $j$.
\item[Transfer:] If the station $j$ has at least one outlet, we consider each station $j', j' \neq j$, and transfer all resources spent on station $j$ in each period $t, t+1, \dotsc, T$ to be spent on station $j'$ instead. To ensure feasibility of the solution, we set the number of outlets at station $j$ for periods $t' \geq t$ to the value in period $t-1$ (or $x_{jk}^0$ if $t=1$). If station $j'$ reaches the maximum number of outlets and there are resources remaining, they are spent on station $j$.
\item[Split:] If the station $j$ has at least one outlet, we consider each station $j', j' > j$, and we evenly split the resources spent on stations $j$ and $j'$ in each period $t, t+1, \dotsc, T$. We note that this move is symmetric, and so it is only necessary to attempt this move if $j' >j$. In order to ensure feasibility of the solution, we can only use this move if the resources spent on these stations and the prior values of the solution for these stations in period $t-1$ allow us to open both stations and place at least one outlet in each.
\end{description}

The moves can be applied to the candidate solution using either the ``first improvement'' or ``best improvement'' methods as described in \cite{Resende2018}. In the first improvement method, the candidate solution $\xhat$ is updated whenever a move is found which improves the objective function $f \left( \xhat\right)$. In the best improvement method, the candidate solution $\xhat$ is updated with the move from amongst all stations $j$ which resulted in the highest objective function $f \left( \xhat\right)$.

\paragraph{Stopping Criteria}

We wish to prevent the local search from spending considerable time making very minor improvements to the candidate solution. At the end of each loop through the stations, we check the increase in the objective function via the local search. If the relative increase to the objective function falls below a given threshold, the local search immediately ends the search in the given period and proceeds to the following one.

The GRASP algorithm terminates when one of the following conditions have been satisfied: i) a threshold number of candidate solutions have been examined,
ii) a threshold number of candidate solutions have been filtered out or,
iii) a time limit has been reached.

\section{Computational Results}
\label{SectionResults}

In this section, we analyse the results from computational experiments using the models and heuristic methods discussed in the previous sections. In Section~\ref{SectionTestEnvironment}, we describe the network and datasets used in our experiments.  {In Section~\ref{SectionSLvsMC}, we assess the advantage of the MC model~\eqref{ModelMaximumCover} over the SL model~\eqref{ModelSinglelevel} when solving both with a MILP solver.} In Section~\ref{subsec:comparison}, we compare the capabilities of the GF model and the MC model. We show that the MC model~\eqref{ModelMaximumCover} can more accurately reflect effects which are known from the literature to affect EV purchase. Finally,  {in Section~\ref{SectionLimitationsMC}, we discuss the limitations of the MC model~\eqref{ModelMaximumCover} motivating,} in Section~\ref{subsec:heuristicsexpriments},  {a computational study comparing it with our} heuristic methods.

All tests were run on a server running Linux version 3.10, with an Intel Core i7-4790 CPU with eight virtual cores and 32 GB of RAM. The code is written in Python 3.7, and is publicly available\footnote{\url{https://github.com/StevenLamontagne/EVChargingStationModel}}. We use CPLEX version 12.10, accessed via the DoCPLEX module (version 2.21). Parameter values can be found in Appendix \ref{AppendixParameterValues}, and instances used in the simulations are publicly accessible\footnote{\url{https://doi.org/10.7488/ds/3850}}.

\subsection{Test Environment}
\label{SectionTestEnvironment}
\begin{figure}
\centering{
{\includegraphics[width=0.5\textwidth]{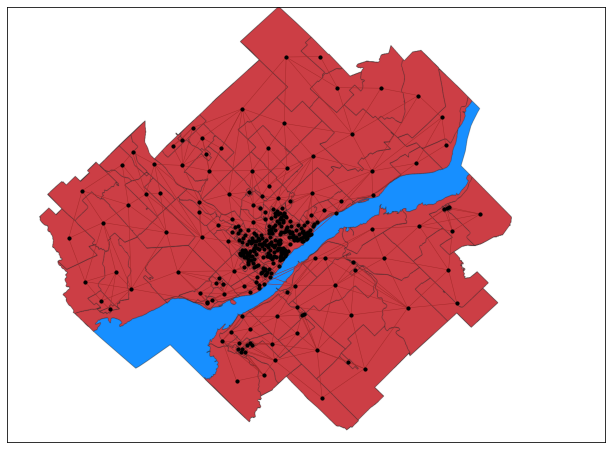}}}
\caption{Trois-Rivières network}
\label{FigureTroisRivieresNetwork}
\end{figure}

The network used for the simulations is based on the smallest aggregation level in the 2016 census~\citep{StatsCan2016} for the city of Trois-Rivières, Qu\'ebec. This defines 317 zones within the city, with populations and aggregate characteristics given for each zone. Nodes in the graph are given by the centroids of each zone. The edges in the graph are created between adjacent zones, with the edge length being the Euclidean distance between the centroids. We note that the Saint Lawrence river divides the city into two parts. Edges have been added to the graph to account for the Laviolette bridge, which connects both parts. The network is shown in Figure~\ref{FigureTroisRivieresNetwork}, with the nodes shown as points. 

We generated five datasets {, each containing 20 instances with} a list of candidate stations, a set of user classes, and pre-computed error terms. The parameters and the generation method for the error terms for each dataset are discussed in Appendix~\ref{AppendixParameterValues}, with key parameter values given in  {Table~\ref{TableKeyParametersMC}}. We give a summary of the important distinctions in each dataset:
\begin{description}
\item[Simple] A small-scale dataset  {where key parameter values are fixed such that the resulting instances are easy to solve.} This is the only dataset which the SL model~\eqref{ModelSinglelevel} can solve. 
\item[Distance] This dataset increases the penalisation term in the  alternative-specific constant for each station to account for distance. 
\item[HomeCharging] For this dataset, we create two user classes in each node: one with and one without access to home charging. 
\item[LongSpan] This dataset increases the number of years from four to ten, but the user classes remain consistent across the planning horizon.
\item[Price] In this dataset, we simulate a decrease in the price of the EV year-by-year, which affects user classes differently based on their income. The alternative-specific constant for each station is modified based on the income level of the user class and the current year.
\end{description}
The Simple, Distance, and HomeCharging datasets have a candidate location set with ten options, and users will not consider any station which is more than ten kilometres away. In contrast, the LongSpan and Price datasets have a candidate location with 30 options, and users may consider any charging station regardless of the distance. This results in significantly more difficult problems to solve. In general, CPLEX is not able to solve the exact models for these instances.

\begin{table}[h]\centering
 {
\resizebox{0.7\textwidth}{!}{\begin{tabular}{l r r r r r}
\hline
Parameter & Simple & Distance & HomeCharging & LongSpan & Price \\
\hline
$T$ & 4 & 4 & 4 & 10 & 4\\
$\abs{M}$ & 10 & 10 & 10 & 30 & 30 \\
$\abs{N}$ & 317 & 317 & 734 & 317 & 1397 \\
$R_i$ & $15-105$ & $15-105$ & $15-105$ & 465 & 465\\
$m_j$ (all stations)& 2  & 6 & 6 & 6 & 6\\
\hline
\end{tabular}
}
\caption{Parameter values for the generated instances. \label{TableKeyParametersMC}}
}
\end{table}

 {
\subsection{Comparison of the Single-Level and Maximum Covering Models}
\label{SectionSLvsMC}
In order to illustrate the advantages of the maximum covering formulation over the single-level version, we solve all instances in the Simple dataset using both formulations. Table~\ref{TableMCvsSL} reports the following statistics averaged over the 20 instances: The CPU time for the root and branch-and-cut, as well as the optimality gaps for each phase and the number of nodes, are as reported in the CPLEX log file. For the number of instances solved, we consider an instance to be solved if the solver terminates with a provably-optimal solution.  For the CPU time in the linear relaxation, neither the Python API nor the CPLEX log directly give the time to solve the linear relaxation of the model. Since CPLEX may run additional methods at the root node to more rapidly solve a discrete optimisation problem (e.g. cuts or heuristics) rather than simply solving its linear relaxation, we ran the solver a second time with the continuous versions of all the variables. As such, the CPU time for the linear relaxation may be higher than for the root node, due to run-time differences. 

We note a large difference in the solving time of the models. In addition to the average improvement between the MC model~\eqref{ModelMaximumCover} and the SL model~\eqref{ModelSinglelevel}, we also note that the former displays small variance: for the SL model~\eqref{ModelSinglelevel}, the solving times ranged from 1,978.21 to 4,724.34 seconds, whereas for the MC model~\eqref{ModelMaximumCover}, the solving times ranged from 0.22 to 0.32 seconds. We also observe improvements in the optimality gap at the root node as well as in the number of nodes explored, likely due to the presence of the Big-M constraints in the SL model~\eqref{ModelSinglelevel}. 

\begin{table}\centering
 {
\resizebox{0.5\textwidth}{!}{\begin{tabular}{lrr}
\hline
                             &  MC &  SL \\
\hline
CPU time (LP, sec) & 0.83 & 4466.75
\\
CPU time (root, sec) & 0.23 & 2294.71
\\
CPU time (B\&C, sec) & 0.03 & 729.24
\\
Optimality gap (root, \%) & 0.69 & 2.15
\\
Optimality gap (B\&C, \%) & 0.00 & 0.00
\\
Nodes explored & 8.50 & 769.85
\\
Instances solved & 20 & 20 
\\
\hline
\end{tabular}
}
\caption{Average performance of CPLEX applied to the MC model and the SL model in the Simple dataset; LP denotes the linear relaxation and B\&C denotes the branch-and-cut. \label{TableMCvsSL}}
}
\end{table}

}%End of  

\subsection{Comparison of the Maximum Covering and GF Models}\label{subsec:comparison}

We compare the capabilities of the GF model~ {\citep{Anjos2020}} with those of the MC model~\eqref{ModelMaximumCover}. To accurately compare the models, we only use the node-based, intracity part of the GF model. Parameter values are chosen to match as closely as possible between both models. The modified GF model, as well as the parameters, are presented in Appendix~\ref{AppendixGrowthFunction}. 

We assume that the capacity of each charging outlet is infinite in the GF model.  This is consistent with the assumption in the MC model that the stations are uncapacitated.  

In Sections~\ref{SectionDistance} and~\ref{SectionHomeCharging}, we consider two cases. In the first case, we force the solver to use the same solution for both the GF and the MC models. This allows for comparing the spread of EVs around charging stations. In the second case, we find the solution returned by the solver for the GF model. We then calculate the objective value of the MC model~\eqref{ModelMaximumCover} using that solution (by using~\eqref{HeuristicSolutionQuality}). This allows us to analyse if the differences in the spread of EVs has an impact on solution quality.

\subsubsection{Distance to Charging Station}
\label{SectionDistance}
Distance is a key factor in determining which facility users choose to patronise in the maximum capture problem~
\citep{Benati2002, Eiselt2019} as well as in models that examine existing EV owners' choice of charging location~\citep{Luo2015, Vermeulen2019, Wolbertus2021}. These works all find that users are less likely to select a facility as the distance increases. To our knowledge, there are no studies which examine the impact of the distance of charging stations to users in the decision to purchase an EV, but we assume that similar results hold in this case. That is, the utility of a charging option (and thus the likelihood that the charging option acts as a primary recharging method) decreases with distance. 

In the GF model, all users have a maximum distance within which they consider charging stations. At a given node, users consider charging at any charging station within that maximum distance. 

In the MC model~\eqref{ModelMaximumCover}, we also assume that the users have a maximum distance for considering charging stations. However, additionally, the utility decreases with distance. For this comparison, we use the Distance dataset described in Section~\ref{SectionTestEnvironment}.

\begin{figure}
\includegraphics[width=0.5\textwidth]{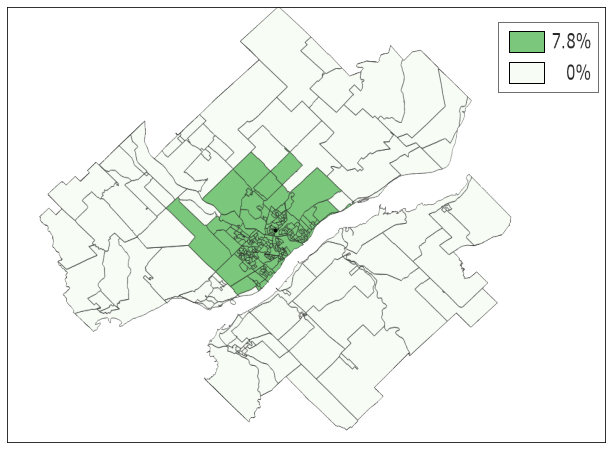}
\includegraphics[width=0.5\textwidth]{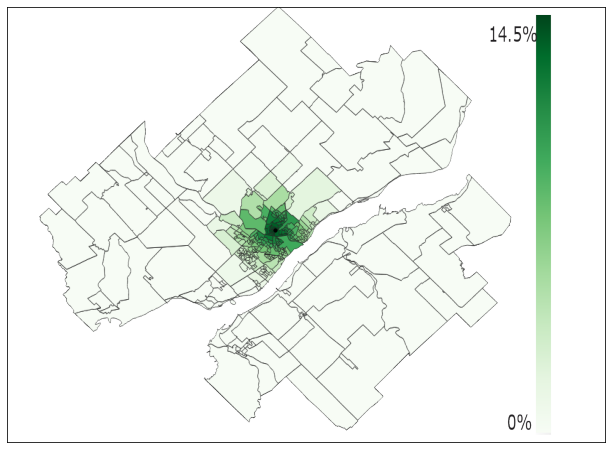}
\caption{Percentage of population in each zone which purchases an EV by the end of the simulation when examining the distance to the charging station. GF model on left side, MC model on right side. \label{FigureComparisonDistance}
}
\end{figure}

In Figure~\ref{FigureComparisonDistance}, we see the percentage of the population that purchases an EV (at the end of the planning horizon) when one station is opened. On the left, in the GF model, we see that the EV adoption rate is the same across the entire region that considers that station. By comparison, on the right in the MC model, we see that the EV adoption rate decreases as the distance increases. 

In Table~\ref{TableObjectivesDistance}, we report the value of~\eqref{HeuristicSolutionQuality} for the solutions from both models. We see that the optimal objective values for the GF model are around 42\% lower than those for the MC model. However, the incentive for the GF model to place more outlets is linked to the capacity of each station. Thus, in the uncapacitated case, it will only place one outlet at any station it opens. To counteract this effect and compare the models more fairly, we examine the objective values of adjusted solutions. These set the capacity of each station selected to be open by the GF model to its maximum capacity of 6 outlets. Despite this being an infeasible solution (due to the budget), the optimal objective values for the adjusted GF model are still around 20.6\% lower than those for the MC model.

\begin{table}\centering\resizebox{0.55\textwidth}{!}{
\begin{tabular}{lrrr}
\hline
{} &  GF &  GF (Adjusted) &  MC \\
\hline
5th percentile  &          9117.18 &                    12954.33 &          16496.45 \\
Median          &          9266.85 &                    13113.74 &          16592.65 \\
95th percentile &          9407.89 &                    13244.96 &          16764.57 \\
\hline
\end{tabular}}
\caption{Number of EVs from the solutions of the GF and MC models.  \label{TableObjectivesDistance}}
\end{table}

We note the difference in the spread of EVs observed in Figure~\ref{FigureComparisonDistance}, with the spread in the MC case being more consistent with the literature  {\citep{Benati2002, Eiselt2019}}. Additionally, as we report in Table~\ref{TableObjectivesDistance}, the quality of the solutions are considerably different between the two models. This indicates that these differences have an important impact in the solutions, and highlights the benefits of the maximum covering formulation.

\subsubsection{Access to Home Charging}
\label{SectionHomeCharging}
Given it is a significant factor in the decision to purchase an EV, accurately modelling access to home charging is critical~\citep{Hidrue2011, Bailey2015}.

%In the GF model, this is done via a constant $a^t$, independent of the zone, which represents the proportion of total EVs that have access to home charging. The total number of users in each zone who may purchase EVs is then bound by the capacity of the stations sufficiently close to them multiplied by $1-a^t$. This has the side effect that if there are no charging outlets sufficiently close to the users, then no EVs may be purchased (not even those that have access to home charging). 

\begin{figure}
\includegraphics[width=0.5\textwidth]{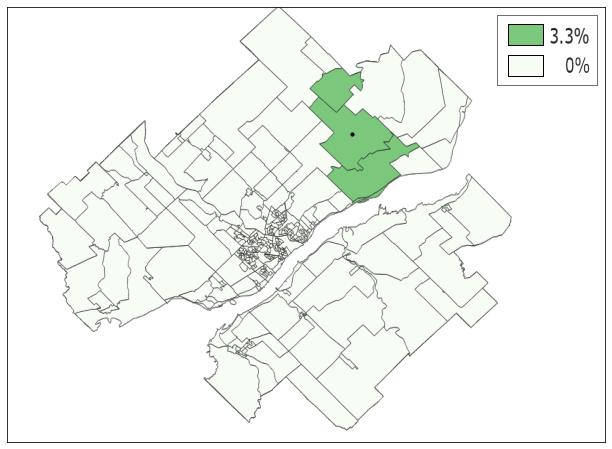}
\includegraphics[width=0.5\textwidth]{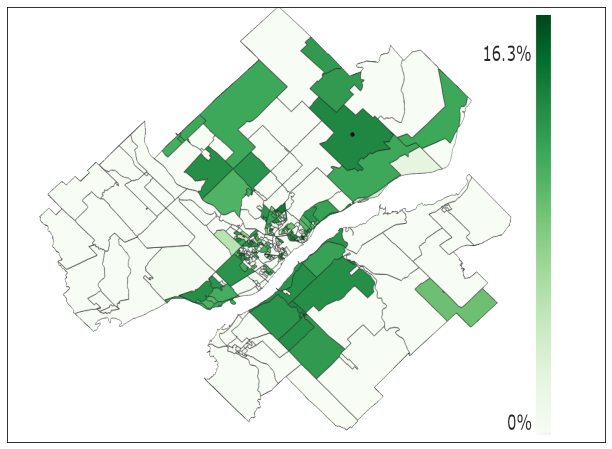}
\caption{Percentage of population in each zone which purchases an EV by the end of the simulation when examining the home charging access. GF model on left side, MC model on right side. \label{FigureComparisonHomeCharging}}
\end{figure}

In the GF model, it is assumed that a given fraction of new EV owners will have access to home charging (which depends on the year, but not the location). The total number of EVs at each location, including both those who have access to home charging and those who do not, is then bound by the capacity of nearby stations. Consequently, if there are no charging outlets sufficiently close to the users, not even those with home charging access may purchase EVs. We note this effect in Figure~\ref{FigureComparisonHomeCharging}.

In the maximum covering model~\eqref{ModelMaximumCover}, the population and home charging access for each user class can be set independently. By creating two user classes for each area, with appropriately set populations, we can more accurately model the percentage of the population that has access to home charging. Additionally, since the two user classes are separate in the model, the users who have access to home charging are allowed to purchase EVs, even if no public charging infrastructure is sufficiently close. For this comparison, we use the HomeCharging dataset described in Section~\ref{SectionTestEnvironment}.

\begin{figure}
\centering
{\includegraphics[width=0.7\textwidth]{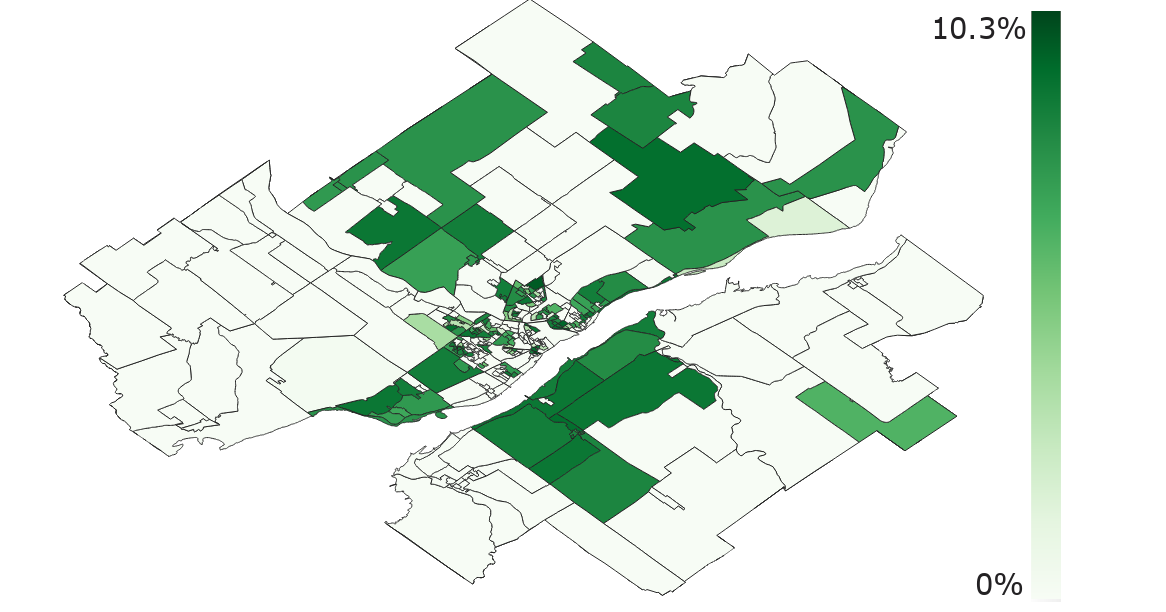}}
\caption{Percentage of Population Covered by Home Charging Access.}
\label{FigureHomeChargingAccess}
\end{figure}

In Figure~\ref{FigureHomeChargingAccess}, we depict the percentage of the population in each area that purchase EVs by the end of the simulation period simply because of home charging access (according to the MC model~\eqref{ModelMaximumCover}). In other words, even if no public charging infrastructure is installed, these users can charge at home and decide to purchase an EV. By contrast, in Figure~\ref{FigureComparisonHomeCharging}, we depict the percentage of the population at the end of the simulation period that purchase EVs when one station is opened. On the left, in the GF model, we note that only users near the charging station have purchased EVs. By comparison on the right, in the MC model~\eqref{ModelMaximumCover}, we note that the users that were covered by home charging in Figure \ref{FigureHomeChargingAccess} have purchased EVs. Additionally, we see that near the  charging station there is an increase to the percentage of the population that purchases EVs, as users who were not covered by home charging access are now covered by the charging station. 

In Table~\ref{FigureObjectivesHomeCharging}, we report the value of~\eqref{HeuristicSolutionQuality} for the solutions from both models. We see that not accounting for home charging access results in a decrease in the number of EVs of around 23.8\%. As before, we examine an adjusted solution to take the infinite capacity into account. Despite this being an infeasible solution, it results in a decrease in the numbers of EVs of around 10.8\% compared to the MC model. Clearly, the same findings emerge as in Section~\ref{SectionDistance}.

\begin{table}\centering\resizebox{0.55\textwidth}{!}{
\begin{tabular}{lrrr}
\hline
{} &  GF &  GF (Adjusted) &  MC \\
\hline
5th percentile  &         13648.43 &                    15990.52 &          17975.55 \\
Median          &         13737.12 &                    16091.57 &          18036.37 \\
95th percentile &         13780.56 &                    16133.63 &          18085.95 \\
\hline
\end{tabular}}
\caption{Number of EVs from the solutions of the GF and MC models.  \label{FigureObjectivesHomeCharging}}
\end{table}

 {
\subsection{Limitations of the Maximum Covering Model}
\label{SectionLimitationsMC}
While the maximum covering formulation significantly improves tractability by MILP solvers compared to the single-level formulation, our problem remains difficult to solve for large instances. To demonstrate this, we solve the MC model~\eqref{ModelMaximumCover} over all instances in all five datasets, restricting CPLEX to a time limit of 7,200 seconds. In Table~\ref{TableMaximumCoveringPerformance}, we report the average performance, using the same statistics reported in Section~\ref{SectionSLvsMC}.

We first note that CPLEX is not able to solve any of the instances in the LongSpan and Price datasets within the time limit. Moreover, not only can it not complete the root node, it is not even able to solve the LP relaxation This highlights the complexity of the problem, even for a moderate-sized city and only 30 candidate locations. Concerning the three smallest datasets (Simple, Distance, and HomeCharging) which are all solved to optimality, we remark that the optimality gap at the root node is quite good, suggesting that the majority of the difficulty lies in solving the root node.  Finally, we note that, while most of the statistics of each instance are close to the average, the number of nodes explored and the CPU time for the branch-and-cut phase vary, depending on how quickly the optimal solution is found. 

\begin{table}\centering
 {
\resizebox{0.8\textwidth}{!}{\begin{tabular}{llrrrrr}
\hline
& Simple & Distance & HomeCharging & LongSpan & Price \\
\hline
CPU time (LP, sec) & 0.83 & 1.01 & 4.33 & -  & -
\\
CPU time (root, sec) & 0.23 & 2.12 & 6.29 & 7202.90 & 7203.88
\\
CPU time (B\&C, sec) & 0.03 & 0.53 & 4.51 & - & -
\\
Optimality gap (root, \%) & 0.69 & 2.47 & 2.41 & 37.34 & 75.14
\\
Optimality gap (B\&C, \%) & 0 & 0 & 0 & 37.34 & 75.14
\\
Nodes explored & 8.50 & 99.95 & 447.75 & 0 & 0 
\\
Instances solved & 20 & 20 & 20 & 0 & 0
\\
\hline
\end{tabular}
}
\caption{Average performance of CPLEX applied to the MC model. The entries for the Simple dataset are copied from Table~\ref{TableMCvsSL}, for ease of comparison. \label{TableMaximumCoveringPerformance}}
}
\end{table}
} %End of  

\subsection{Comparing Heuristics}
\label{subsec:heuristicsexpriments}

We compare solving the MC model~\eqref{ModelMaximumCover} via the standard branch-and-cut approach used by CPLEX (labelled as ``Exact'' in the results) and the heuristics presented in Section~\ref{SectionHeuristics}.  {We impose a time limit of two hours (7,200 seconds) for the Exact and all heuristic methods. In the cases of the Exact method and the Rolling Horizon method, this time limit is given directly to CPLEX, and thus only limits the CPLEX solving time. } The Greedy procedure is tested using both the myopic and hyperoptic search methods. For the Rolling Horizon method, we examine two ways of distributing the time across the years. In the first, we divide the time evenly amongst all years (``Even''). In the second, the time limit is divided geometrically, with the first year having a time limit of 3,600 seconds and each successive year having half of the previous time (``Geometric'').  For the GRASP procedure, we test both the myopic and hyperoptic search methods.  We set the value of $\alpha$ to 0.85, as recommended for the maximum covering problem in~\cite{Resende2018}. The GRASP procedures run until 300 solutions have been examined, or until 500 solutions have been filtered out. The second phase of the GRASP procedure uses the first improvement method.
  
 %end of  
\subsubsection{Solving Time}

In Table~\ref{TableHeuristicTimes}, we report the 5th percentile, the 95th percentile, and the average solving time across all instances of a given dataset and method. For the small datasets (Simple, Distance, HomeCharging), the solving time of the Exact method is only surpassed by the Greedy  {and Rolling Horizon} methods. Both Greedy methods are undoubtedly the fastest. This is the case in the Simple and Distance datasets  {(with the Rolling Horizon methods being nearly as fast in the Simple dataset), however, it is more noticeable in the HomeCharging, Price, and LongSpan datasets. }

The Rolling Horizon methods both solve quickly for the small datasets. In the LongSpan test set, the Rolling Horizon methods solve more quickly than the Exact method. Since the user classes are the same for each year, solutions from previous years act as good quality warmstart solutions. This allows the later years in the Rolling Horizon method to terminate before the time limit. However, this does not occur in the Price test set. Since the user classes change each year, the quality of the solution for previous years is poorer, which causes the solver to reach the time limit. 

While the GRASP methods took longer to solve than the RollingHorizon methods in the small datasets, they scaled better for the Price dataset. This resulted in a solving time around a sixth that of the RollingHorizon and Exact methods. In the LongSpan dataset, the solving times for the GRASP method are relatively high. This is due to the fact that both the first and second phases iterate over every year, which causes the solving time to increase substantially as the number of time years increases. The hyperoptic method is solved noticeably slower than the myopic method.

\begin{table}
 {
\centering
\resizebox{0.8\textwidth}{!}{\begin{tabular}{llrrrrr}
\hline
                   &                 &  Simple &  Distance &  HomeCharging &   Price &  LongSpan \\
\hline
Exact & 5th percentile &    0.22 &      2.26 &          9.17 & 7203.50 &   7202.70 \\
                   & Average &    0.27 &      2.66 &         10.80 & 7203.98 &   7202.99 \\
                   & 95th percentile &    0.30 &      3.04 &         11.89 & 7204.51 &   7203.27 \\
\hline
GreedyOne (Myopic) & 5th percentile &    0.10 &      0.15 &          0.15 &    1.03 &      1.17 \\
                   & Average &    0.10 &      0.16 &          0.39 &    1.06 &      1.23 \\
                   & 95th percentile &    0.10 &      0.17 &          1.07 &    1.09 &      1.65 \\
\hline
GreedyOne (Hyperoptic) & 5th percentile &    0.21 &      0.29 &          0.29 &    2.20 &      4.98 \\
                   & Average &    0.21 &      0.31 &          0.57 &    2.27 &      5.09 \\
                   & 95th percentile &    0.21 &      0.33 &          1.19 &    2.36 &      5.12 \\
\hline
RollingHorizon (Even) & 5th percentile &    0.17 &      0.37 &          1.92 & 7201.49 &   3557.43 \\
                   & Average &    0.18 &      0.40 &          2.09 & 7226.27 &   3919.75 \\
                   & 95th percentile &    0.19 &      0.41 &          2.22 & 7298.76 &   4238.00 \\
\hline
RollingHorizon (Geometric) & 5th percentile &    0.17 &      0.37 &          1.92 & 7201.47 &   7096.51 \\
                   & Average &    0.18 &      0.40 &          2.08 & 7219.97 &   7163.96 \\
                   & 95th percentile &    0.19 &      0.41 &          2.21 & 7314.78 &   7186.91 \\
\hline
GRASP (Myopic) & 5th percentile &   66.48 &     78.20 &         88.06 &  851.76 &   1593.75 \\
                   & Average &   71.22 &     83.74 &         91.55 &  873.47 &   1635.97 \\
                   & 95th percentile &   81.70 &     87.18 &        101.45 &  896.95 &   1666.46 \\
\hline
GRASP (Hyperoptic) & 5th percentile &   99.27 &    119.29 &        128.16 & 1405.71 &   3058.92 \\
                   & Average &  105.07 &    122.85 &        132.76 & 1433.38 &   3121.60 \\
                   & 95th percentile &  113.94 &    127.47 &        138.17 & 1478.21 &   3197.92 \\
\hline
\end{tabular}
}
\caption{Comparison of solving times for heuristic methods across all datasets. Times are given in seconds. \label{TableHeuristicTimes}}
} %End of  
\end{table}

\subsubsection{Solution Quality}

In Table~\ref{TableHeuristicGaps}, we report the 5th percentile, the 95th percentile, and the average gaps to the best known solution for all instances of a given test set and method.  { More specifically, for each instance, we examine which method found the solution with the highest objective, whose value is denoted as $bestSol$. The gap for each method is then given by
$
gap(x^\star) = \frac{bestSol - f(x^\star)}{bestSol},
$
where $x^\star$ is the solution found by the given method. 
}

Additionally, we give the number of instances in which a method has found a solution with objective value equal to that of the best known solution.  In the case of the Simple, Distance, and HomeCharging datasets, the exact solution is known, thus the corresponding entries for the heuristic methods are the optimality gaps.

While the Exact method produces the best solution value for the smaller tests (Simple, Distance, HomeCharging),  {it is unable to find solutions which are as good as the heuristic methods in the Price and LongSpan datasets.}  

The Greedy methods performed slightly worse than the other heuristic methods in the Simple, Distance, and Home 
Charging datasets, but they performed surprisingly well in the LongSpan and Price datasets. They frequently found the best solution, and with gaps under 0.1\% even in cases where the method did not find the best solution. 

The Rolling Horizon methods perform well in the smaller test sets but, similar to the Exact method,  {they do not perform as well in the LongSpan and Price datasets as the other heuristic methods. }  
Both GRASP methods performed comparable to each other, and comparable with the Greedy methods.  {However, the GRASP methods were more consistently able to produce good-quality solutions across all datasets as compared to the Greedy methods, as can be seen in the Distance dataset. }

\begin{table}
 {
\centering
\resizebox{0.8\textwidth}{!}{\begin{tabular}{llrrrrr}
\hline
                   &            &  Simple &  Distance &  HomeCharging &  Price &  LongSpan \\
\hline
Exact & 5th percentile &    0.00 &      0.00 &          0.00 &  12.67 &      8.86 \\
                   & Average &    0.00 &      0.00 &          0.00 &  12.74 &      8.94 \\
                   & 95th percentile &    0.00 &      0.00 &          0.01 &  12.83 &      9.01 \\
                   &  \# of best &   20 &     20 &         20 &   0 &      0 \\
\hline
GreedyOne (Myopic) & 5th percentile &    0.00 &      1.66 &          0.46 &   0.00 &      0.00 \\
                   & Average &    0.02 &      2.08 &          0.53 &   0.01 &      0.01 \\
                   & 95th percentile &    0.07 &      2.76 &          0.68 &   0.05 &      0.03 \\
                   &  \# of best &   14 &      0 &          0 &  11 &     15 \\
\hline
GreedyOne (Hyperoptic) & 5th percentile &    0.00 &      1.47 &          0.44 &   0.00 &      0.00 \\
                   & Average &    0.07 &      2.03 &          0.61 &   0.02 &      0.02 \\
                   & 95th percentile &    0.19 &      2.53 &          0.84 &   0.06 &      0.05 \\
                   &  \# of best &    9 &      0 &          0 &   8 &     10 \\
\hline
RollingHorizon (Even) & 5th percentile &    0.00 &      0.00 &          0.00 &  12.33 &      4.11 \\
                   & Average &    0.00 &      0.07 &          0.02 &  12.66 &      4.19 \\
                   & 95th percentile &    0.00 &      0.38 &          0.06 &  12.83 &      4.27 \\
                   &  \# of best &   20 &     14 &         12 &   0 &      0 \\
\hline
RollingHorizon (Geometric) & 5th percentile &    0.00 &      0.00 &          0.00 &  12.51 &      2.52 \\
                   & Average &    0.00 &      0.07 &          0.02 &  12.72 &      3.95 \\
                   & 95th percentile &    0.00 &      0.38 &          0.06 &  12.83 &      4.26 \\
                   &  \# of best &   20 &     14 &         12 &   0 &      0 \\
\hline
GRASP (Myopic) & 5th percentile &    0.00 &      0.14 &          0.16 &   0.00 &      0.10 \\
                   & Average &    0.02 &      0.43 &          0.26 &   0.02 &      0.18 \\
                   & 95th percentile &    0.12 &      0.68 &          0.36 &   0.04 &      0.24 \\
                   &  \# of best &   16 &      0 &          0 &   5 &      0 \\
\hline
GRASP (Hyperoptic) & 5th percentile &    0.00 &      0.19 &          0.12 &   0.00 &      0.10 \\
                   & Average &    0.03 &      0.44 &          0.23 &   0.04 &      0.19 \\
                   & 95th percentile &    0.16 &      0.67 &          0.36 &   0.07 &      0.26 \\
                   &  \# of best &   14 &      0 &          0 &   3 &      0 \\
\hline
\end{tabular}
}
\caption{Comparison of gaps to the best known solution for heuristic methods in all five datasets (in percentage). \label{TableHeuristicGaps}}
} %End of  
\end{table}

\section{Conclusion}
\label{SectionConclusion}

In this work, we proposed a model for determining the optimal location of EV charging stations in a long-term planning environment so as to maximise the total number of EVs. To consider user-specific characteristics, we used discrete choice models to represent the decision of the users to purchase EVs.  When compared to the existing model for this problem \citep[the GF model in][]{Anjos2020} this allowed for more intuitive user behaviour with regards to charging station location and home charging access. Additionally, solutions for the GF model were significantly different than those for the MC model, highlighting the benefits of using the latter.

Using the simulation-based approach of~\cite{Paneque2020} resulted in a bilevel model. By reformulating the model as a maximum covering problem, the sets of Big-M constraints were removed, thus improving the tractability. We note that this reformulation can be applied to other applications of the simulation-based method (including the MCPRU) if no capacity constraints are present. 

For more difficult instances, several heuristic methods were proposed for obtaining feasible solutions. While the Rolling Horizon heuristic was not as effective for more difficult instances, the Greedy and GRASP algorithms both performed very well. In particular, we note that the Greedy methods were able to find solutions within a few seconds, even on the most difficult instances.  {Additionally, the GRASP methods can obtain solutions of  near-optimal quality within a few minutes for the small instances, and they can easily surpass the quality of the best solution obtained from the exact method, within less than half the time limit.} The similar performances of the Greedy and GRASP methods suggests that more complex local search moves may be necessary to find better solutions.

While the heuristic methods were able to quickly find feasible solutions, we  are not able to verify the quality of these solutions on the larger instances.  Further research will examine an exact method specifically tailored for larger scale instances.  {A few approaches are possible. Firstly, the structure of the problem resembles the resource-constrained production scheduling problem (RCPSP), for which a specialised algorithm exists~\citep{Bienstock2010, Munoz2018}. The method (which is similar to column generation) has a mechanism for bounding the number of basis elements. However, when this bound is applied to our problem, it is greater than the maximum number of basis elements. As such, a complete enumeration is possible, and the method may require solving the full MC model~\eqref{ModelMaximumCover}. A successful modification and implementation of this algorithm for our problem would be the subject of future research. A second possible approach is the Bender's decomposition method proposed in~\cite{Cordeau2019}, which is designed for large-scale maximum covering problems. Since the restricted master problem does not include the full set of covering constraints, and solutions to the subproblem can be calculated analytically, the method seems quite well-suited. However, in our preliminary testing, we have found poor convergence. In part, this is likely due to the higher ratio of stations to users as compared with the original paper, which the authors note to be an important factor. As such, research is ongoing to improve the performance of this method in our problem.
}

Since the MC model~\eqref{ModelMaximumCover} was designed for the intracity network, it is not applicable to larger geographical areas. For example, in the MC model~\eqref{ModelMaximumCover}, there are no constraints related to the range of EVs, which would be necessary for intercity travel. This is a notable difference with the GF model of \cite{Anjos2020}, which was designed for intracity and intercity networks, and was applied to the province of Qu\'ebec. As such, further research will examine an extension of the MC model~\eqref{ModelMaximumCover} which includes the intercity network.

% Acknowledgments here
\section*{Acknowledgements}

 The authors thank the anonymous reviewers whose comments helped to improve the manuscript. In particular, one reviewer proposed a reformulation of the model which greatly improved the convergence rate, and led to the version presented here.

The authors gratefully acknowledge the assistance of Jean-Luc Dupr\'e from \emph{Direction Mobilit\'e} of \emph{Hydro-Qu\'ebec}, for both providing and helping to  understand EV charging data, as well as for sharing his expertise on EV charging stations and the network. We gratefully acknowledge the work of Mahsa Moghaddass in estimating the parameter values used in the models.

This research was supported by Hydro-Québec, NSERC Collaborative Research and Development Grant CRDPJ 536757 - 19, and NSERC Discovery grant 2017-06054.

\bibliographystyle{plainnat} % outcomment this and next line in Case 1
\bibliography{References/Benders,References/Bilevel, References/ChargingStationLocation, References/ChargingStationChoice, References/Data, References/DiscreteChoiceModels, References/EVAcquisition, References/MaximumCapture, References/MaximumCovering, References/Heuristics, References/Review} % if more than one, comma separated

% Appendix here
% Options are (1) APPENDIX (with or without general title) or 
%             (2) APPENDICES (if it has more than one unrelated sections)
% Outcomment the appropriate case if necessary
%
% \begin{APPENDIX}{<Title of the Appendix>}
% \end{APPENDIX}
%
%   or 
%
\clearpage
\appendix

\section{Bilevel Optimisation Model}
\label{AppendixBilevelModel}
The natural hierarchical structure of our problem suggests a bilevel formulation. Bilevel optimisation models sequential decision making, where first, the leader  takes a decision (\emph{upper level}) and then, the followers react by solving an optimisation problem (\emph{lower level}). The optimal solution of a bilevel model are the decision values for the leader that optimise its objective function, based on the optimal reaction of the lower-level users to the values of those variables.  {In our case, the decision maker is in the upper-level while the users are in the lower-level.} 

The choice in the lower-level to purchase an EV is modelled via a discrete choice model, where
 {the users (followers) maximise their \emph{utility}.} More specifically, the users must choose an alternative from a finite set of available alternatives (here, open stations, home charging, and opt-out). The value of each alternative is predicted through the use of a \emph{utility function}, which associates the value of a given alternative for users based on observable and unobservable factors. 
Under the RUM assumption, users then, as rational beings, select the alternative which presents the maximum benefit to them, as represented by the alternative with the highest utility. For each period $1 \leq t \leq T$, user class $i \in N$, and alternative $j \in  {\mathcal{C}_i^t(\x)}$, the utility is denoted $u_{ji}^t$.

The analyst has imperfect knowledge of the utility of the users, so we model it as a random variable. Hence, instead of a deterministic model identifying the alternative chosen by the users, we obtain a probability distribution over the set of available alternatives. Consequently, the leader maximises the expected number of users purchasing an EV or, equivalently, minimises the expected number of users that \emph{do not} purchase an EV. This is given by 
\begin{equation}
\label{BilevelGeneralProblemProbability}
\min_{x \in X} \sum_{t =1}^T \sum_{i \in N} N_i^t \mathbb{P}_\epsilon \left[ u_{0i}^t (x, \epsilon ) \geq u_{ji}^t (x, \epsilon ), \forall j \in  {\mathcal{C}_i^t (x)}\right],
\end{equation}
where $X$ denotes the upper-level constraints, and $\epsilon$ denotes the random error term. These constraints are discussed in more detail in Appendix~\ref{SectionBilevelUpper}. We recall that index $j=0$ indicates the opt-out alternative, thus $u_{0i}^t$ is the opt-out utility for user class $i$ in year $t$. 

If the error terms are independent, and identically extreme-value type I distributed, the choice probabilities have an analytic formula (the well-known MNL model). In our application, this assumption may not hold. For example, if stations are near each other, they will also have similar amenities near them (e.g. restaurants, shopping centres, etc.). If a user places high value in those amenities (and they are not explicitly included in the observable factors), the error terms for those stations may be highly correlated. To allow for general discrete choice models which relax this restriction and support flexible substitution patterns, we use the simulation-based approach of \cite{Paneque2020}.

 {Recall that} $R_i$ is the number of \emph{scenarios} for user class $i$,  {$P$ is the set of triplets $(t,i,r)$ for user class $i$, alternative $j$, scenario $r$, and period $t$}, and $\epsilon_{ji}^{rt}$  {is} the realisation of the random variable $\epsilon$ in triplet $(t,i,r)$ . Then, $\forall j \in  {C_i^t(x)}$,  we denote 
\begin{equation}
\label{BilevelWVariables}
w_{ji}^{rt} = \left\{
\begin{aligned}
&1, \text{ if } j \in \argmax_{j' \in  {C_i^t(x)}} \brace{u_{j'i}^t (x, \epsilon_{j'i}^{rt})},
\\
&0, \text{ otherwise.}
\end{aligned}
\right.
\end{equation}
The details of the lower-level problem are discussed in Appendix~\ref{SectionBilevelLower}. For the sake of simplicity, in a flagrant abuse of notation, we denote the vector $\w = \brace{w_{ji}^{rt}}$ with $w_{ji}^{rt}$ given by~\eqref{BilevelWVariables} for each $(t,i,r), j \in  {C_i^t(x)}$ as $\w \in \argmax_{ {j \in C(x)}} \brace{u(x, \epsilon)}$. Then, we can write a sample average approximation of~\eqref{BilevelGeneralProblemProbability} as
\begin{align*}
\min_{x \in X} & \sum_{(t,i,r) \in P} \frac{N_i^t}{R_i}  w_{0i}^{rt},  
\\
\text{s.t. }& \w = \argmax_{j \in  {C(x)}} \brace{u(x, \epsilon)}. 
\end{align*}

\subsection{Lower-Level Problem}
\label{SectionBilevelLower}
The lower-level problem is assumed to be separable for each user class, each period, and each scenario, meaning that there is no interaction among them. Therefore, in what follows, we concentrate on detailing a given triplet $(t,i,r)$. For each alternative $j \in \mathcal{C}_i^t (x)$, let $u_{ji}^{rt} = u_{ji}^t (x, \epsilon_{ji}^{rt})$ be the \emph{simulated utility}.

We previously remarked that the choice set $\mathcal{C}_i^t(\x)$ depends on which stations are open  {and we defined the  choice sets $C_i^{0t}$ and $C_i^{1t}$ as the sets related to alternatives exogenous and endogenous to the optimisation model, respectively}.  
To ensure that the alternative associated with a closed station $j \in C_i^{1t}$ cannot be chosen, we set the simulated utility $u_{ji}^{rt}$ to a lower bound if  {$x_{j1}^t =0$}. This concept is referred to as the ``discounted utility'' in~\cite{Paneque2020}. Let  {$\underline{a}_i^{t}= \min  \left( \brace{ \kappa_{ji}^t +\epsilon_{ji}^{rt}, j \in C_i^{1t}, 1 \leq r \leq R_i} \right)$} and $b_{ji}^{rt} =\sum_{k=1}^{m_j} \beta_{jik}^t +\kappa_{ji}^t +\epsilon_{ji}^{rt}$ be respectively lower and upper bounds on the simulated utility $u_{ji}^{rt}$, and  
 let $\nu_{ji}^{t}=b_{ji}^{rt}-\underline{a}_{i}^{t}$. For each $i \in N$, we assume that $R_i$ is sufficiently large such that for each $1 \leq t \leq T, 1 \leq r \leq R_i$ we have $\underline{a}_i^{t} < u_{0i}^{rt}$.  {We note that the lower bound $\underline{a}_i^{t}$ could be strengthened by adding the utility of one outlet, (e.g. $\min \brace{\beta_{j1}^{rt} + \kappa_{ji}^t +\epsilon_{ji}^{rt}}$). However, in our testing, this did not have a significant impact. We also note that, since the error terms $\epsilon_{ji}^{rt}$ can come from unbounded distributions, it is not generally possible to use a fixed lower bound for $u_{0i}^{rt}$, e.g. $\underline{a}_{i}^{t}=0$.}

The linear formulation for the simulated utility $u_{ji}^{rt}$ is given by
\begin{align}
u_{ji}^{rt} & \geq  {\underline{a}_i^t},  && j \in C_i^{1t}, (t,i,r) \in P,  \label{ConstraintDiscountClosed1Bilevel}
\\
u_{ji}^{rt} & \leq  {\underline{a}_i^t} + \nu_{ji}^{rt} x_{j1}^t, &&  j \in C_i^{1t}, (t,i,r) \in P,  \label{ConstraintDiscountClosed2Bilevel}
\\
u_{ji}^{rt} & \geq \sum_{k=1}^{m_j} \beta_{jik}^t x_{jk}^t+\kappa_{ji}^t +\epsilon_{ji}^{rt} - \nu_{ji} ^{rt} \left( 1 - x_{j1}^t \right), &&  j \in C_i^{1t}, (t,i,r) \in P,  \label{ConstraintDiscountOpen1Bilevel}
\\
u_{ji}^{rt} & \leq \sum_{k=1}^{m_j} \beta_{jik}^t x_{jk}^t+\kappa_{ji}^t +\epsilon_{ji}^{rt}, &&  j \in C_i^{1t}, (t,i,r) \in P.  \label{ConstraintDiscountOpen2Bilevel}
\end{align}

For each  {$(t,i,r) \in P$}, the value of $w_{ji}^{rt}$ for $j \in C_i^{0t} \cup C_i^{1t}$ is then given by the solution of the following optimisation problem, which acts as the lower-level problem in our bilevel optimisation model:
\begin{subequations}
	\begin{alignat}{5}
\operatorname{Maximise} &\sum_{j \in C_i^{0t}} w_{ji}^{rt}u_{ji}^{rt}+\sum_{j \in C_i^{1t}} w_{ji}^{rt}u_{ji}^{rt}, \\
\text{subject to }&\sum_{j \in C_i^{0t}} w_{ji}^{rt}+\sum_{j \in C_i^{1t}}w_{ji}^{rt}=1, \label{Loweleve:const1}
\\
& w_{ji}^{rt} \in \brace{0,1}, \quad j \in C_i^{0t} \cup C_i^{1t}. \label{Loweleve:const2}
	\end{alignat}
\label{ConstraintLowerLevelProblemBiLevel}
\end{subequations}
It is easy to see that the binary requirements can be relaxed.

\subsection{Upper-Level Problem}
\label{SectionBilevelUpper}

The placement of charging outlets and stations in the upper-level is restricted by the following set of constraints
 {
\begin{align}
&\sum_{j \in M}\sum_{k=1}^{m_j} c_{jk}^ t\left(x_{jk}^ t-x_{jk}^{t-1} \right)  \leq B^t, && 1 \leq t \leq T. \label{ConstraintBudgetBiLevel} 
\\
&x_{jk}^t \leq x_{jk-1}^t,  && 1 \leq t \leq T, j \in M, 
\label{ConstraintPayOneTimeCostBiLevel} 
\\
&x_{jk}^t \geq x_{jk}^{t-1}, && 1 \leq t \leq T, j \in M, 1 \leq k \leq m_j. 
\label{ConstraintCantCloseStationsBiLevel}
\end{align}
These correspond to Constraints~\eqref{ModelSinglelevel:Budget}-\eqref{ModelSinglelevel:KeepOutlets} presented in Section~\ref{subsec:singlelevel}.}
%Constraints~\eqref{ConstraintBudgetBiLevel} ensure that the amount spent installing charging outlets does not exceed the budget for that period. 
Note that it would also be possible to supplement the per-period budget with an overall budget, as was done in \cite{Anjos2020}.
% end of  

 {Constraints~\eqref{ConstraintPayOneTimeCostBiLevel} enforces that if we have at least $k$ outlets, we must also have at least $k-1$ outlets.}

Constraints~\eqref{ConstraintCantCloseStationsBiLevel} forbid the model from removing charging outlets. These constraints assume that it would be suboptimal to remove a station. 

\subsection{Bilevel Model}
\label{SectionBilevelBilevelModel}
We now introduce the full, bilevel model
\begin{align}
\operatorname{Minimise} & \sum_{(t,i,r) \in P}  \frac{N_i^{t}}{R_i} w_{0i}^{rt}, \label{Obj:bilevel}
\\
\text{subject to } &  \eqref{UtilityCalc_Exo}, \eqref{ConstraintDiscountClosed1Bilevel}-\eqref{ConstraintDiscountOpen2Bilevel}, \eqref{ConstraintBudgetBiLevel} -\eqref{ConstraintCantCloseStationsBiLevel} \notag \\ 
& w_{ji}^{rt} \in \argmax\Bigg\{  \sum_{j \in C_i^{0t}} w_{ji}^{rt}u_{ji}^{rt}+\sum_{j \in C_i^{1t}} w_{ji}^{rt}u_{ji}^{rt}:   \eqref{Loweleve:const1}-\eqref{Loweleve:const2} \Bigg\} \nonumber
\\
&  {u_{ji}^{rt} \in \R,} \nonumber
\\
& x_{jk}^t\in \brace{0, 1}. \nonumber
\end{align}

 We consider the \emph{optimistic} version of the bilevel problem which means that the users do not select the opt-out alternative if a different alternative has equal utility. While this in theory has a zero probability (given that the error terms are drawn from continuous distributions), this can occur in practice due to numerical precision.

In order to solve the model, we reformulate it as a single-level optimisation problem by transforming the lower-level model~\eqref{ConstraintLowerLevelProblemBiLevel} into a series of constraints for the upper-level model. To this end, we apply the Karush-Kuhn-Tucker conditions which are necessary and sufficient for the optimality of the (linear) lower-level problem~\citep{Sinha2017}, and we linearize the terms $w^{rt}_{ji} \cdot u^{rt}_{ji}$ through Big-M constraints. In this way, for each $(t,i,r) \in P$,  the lower-level problem~\eqref{ConstraintLowerLevelProblemBiLevel} is replaced by the following constraints:
\begin{align}
& u_{ji}^{rt}-\alpha_{i}^{rt}+\left(1-w_{ji}^{rt}\right)\mu_{ji}^{rt}  \geq 0, \, j \in C_i^{0t} \cup C_i^{1t}, \label{LowerLevelSolutionStartSubproblem}\\
& \sum_{j \in C_i^{0t}} w_{ji}^{rt}+\sum_{j \in C_i^{1t}} w_{ji}^{rt}=1, \\
& \alpha_{i}^{rt} \geq u_{ji}^{rt}, \, j \in C_i^{0t} \cup C_i^{1t}, \\
& w_{ji}^{rt} \in \{0,1\}, j \in C_i^{0t} \cup C_i^{1t}, \\
& \alpha_{i}^{rt} \in \mathbb{R}, \label{LowerLevelSolutionEndSubproblem}
\end{align}
where the Big-M constants, $\mu_{ji}^{rt}$, are given by
\begin{align*}
\mu_{ji}^{rt} &=
\Bigg\{ 
\begin{aligned}
&\max \left( \brace{ b_{ji}^{rt}, j \in C_i^{1t}} \cup \brace{ \kappa_{ji}^{rt} + \epsilon_{ji}^{rt}, j \in C_i^{0t} } \right) - \kappa_{ji}^{rt} - \epsilon_{ji}^{rt}, \quad j \in C_i^{0t}, \\
&\max \left( \brace{ b_{ji}^{rt}, j \in C_i^{1t
}} \cup \brace{ \kappa_{ji}^{rt} + \epsilon_{ji}^{rt}, j \in C_i^{0t} } \right) -  {\underline{a}_i^{t}}, \quad j \in C_i^{1t}.\\
\end{aligned}
\end{align*}
%We refer to the complete model with the constraints~\eqref{LowerLevelSolutionStartSubproblem}-\eqref{LowerLevelSolutionEndSubproblem} as the `Single-Level model' in order to distinguish from the general context of single-level reductions. 

\section{Parameter Values}
\label{AppendixParameterValues}
 
In this section, we describe the parameter values for each dataset, as well as the mechanism for drawing error terms  {for each of the instances}. We start by describing the general framework used for the error terms, as that is common to all of the datasets. Table~\ref{TableParametersMC} provides a list of parameter values, with more detailed explanations in the following subsections for parameter values which differ. Unless otherwise specified, parameter values are set arbitrarily.  { We note that, in all datasets, the number of scenarios is set to $15 \times \abs{C_0^{it} \cup C_1^{it}}$. In the case of the Simple, Distance, and HomeCharging datasets, the size of $C_1^{it}$ varies by user class, due to the maximum distance of 10km for considering a charging station. As a consequence, the number of scenarios varies from 15 (only the opt-out is considered) to 105 (opt-out plus six charging stations). On the contrary, in the Price and LongSpan datasets, there is no maximum distance for consideration, and thus every station is included in the set $C_1^{it}$. As a consequence, $\abs{C_0^{it} \cup C_1^{it}} = 31$, and the number of scenarios is always $465$.

}

In all of our datasets:
\begin{itemize}
\item Each user class $i \in N$ includes the home location (as a node in the network). This allows us to estimate population based on the  census data~\citep{StatsCan2016}, with the number of user classes per node and the partitioning method depending on the dataset. 
\item  {
Each dataset includes 20 instances, where each instance generates different sets of error terms $\epsilon_{ji}^{rt}, j \in C_{0i}^{rt} \cup C_{1i}^{rt}$.
}
\item For each $i \in N, 1 \leq t \leq T$, the alternative-specific constant for the opt-out option $\kappa_{0i}^{t}$ is set to 4.5. 
\end{itemize}

In order to simulate the error terms for the demand model, we employ the error components formulation of the mixed logit model to approximate a nested logit model, as described in  \cite{Train2002} and \cite{Walker2004}. The notation in what follows matches the latter work, and we refer to the aforementioned work for detailed explanations of the process.

For each $i \in N, 1 \leq t \leq T, 1 \leq r \leq R_i$, the vector of error terms $\epsilon_{i}^{rt}=(\epsilon_{i}^{rt})_{j \in C_i^{0t} \cup C_i^{1t}}$ is given by
\begin{equation}
\epsilon_{i}^{rt}= F T \xi^{r} + \zeta^{r},
\end{equation}
with
\begin{itemize}
\item $F$ a \emph{factor loading} matrix.
\item $T$ a diagonal matrix with the standard deviation of each factor.
\item $\xi^{r}$ a vector of IID random terms from a normal distribution. 
\item $\zeta^{r}$ a vector of IID random terms from a Gumbel distribution. 
\end{itemize}
The form of the matrices $F$ and $T$ vary in each dataset. However, in all datasets, $\xi^{r}$ has a location of zero and a scale of one, and $\zeta^{r}$ has location of zero and a scale of three.

\begin{table}[h]\centering
\resizebox{0.7\textwidth}{!}{\begin{tabular}{l r r r r r}
\hline
Parameter & Simple & Distance & HomeCharging & LongSpan & Price \\
\hline
$T$ & 4 & 4 & 4 & 10 & 4\\
$\abs{M}$ & 10 & 10 & 10 & 30 & 30 \\
$\abs{N}$ & 317 & 317 & 634 & 317 & 1397 \\
$\abs{C_0^{it}}$ & 1 & 2 & 1 & 1 & 1 \\
$\abs{C_1^{it}}$ & Varies & Varies & Varies & 30 & 30 \\
$x_{jk}^0$ (all stations and outlets)& 0 & 0 & 0 & 0 & 0  \\
$R_i$ &  {$15-105$} &  {$15-105$} &  {$15-105$} & 465 & 465\\
$B^t$ (per year)& 400 & 400 & 400 & 400 & 400 \\
$m_j$ (all stations)& 2  & 6 & 6 & 6 & 6\\
$c_{j1}^t$ (all stations and years)& 150 & 150 & 150 & 150 & 150\\
$c_{jk}^t$ (all other outlets)& 50  & 50 & 50 & 50 & 50\\
\hline
\end{tabular}
}
\caption{Parameter values for the generated instances \label{TableParametersMC}}
\end{table}

\subsection{Simple Dataset}
The set of user classes $N$ includes one user class for every node in the network. The population of user class $i$, $N_i^t$, is given by the population of the node in the 2016 census multiplied by a factor of 0.1. In other words, 10\% of the population in each node are deciding to purchase a vehicle each year.

For each $i \in N, 1 \leq t \leq T$, the choice set $C_i^{1t}$ includes all stations which are within ten kilometres of the location of the user class.
The utility for $j \in C_i^{1t}$ is linear in terms of the number of charging outlets, with 
\begin{equation}
\beta_{jik}^{t} = 0.281 k, \forall j \in C_i^{1t}.
\end{equation}
Additionally, the alternative-specific constant for each station $j \in C_i^{1t}$ is calculated as
\begin{equation}
\kappa_{ji}^{t} = 1.464 \delta_1 - 0.063 \delta_2 + 0.174 \delta_3,
\end{equation}
with
\begin{itemize}
\item $\delta_1$: binary coefficient indicating if the station is level 3 (i.e. fast charging). We note that in our tests all stations were considered level 3.
\item $\delta_2$: the distance (in kilometres, shortest path in the network) between the user's home and the charging station, 
\item $\delta_3$: binary coefficient indicating if the station is in the city center (defined as a subset of the nodes in the network).
\end{itemize} 
The coefficients for these parameters were estimated using real-world data. A discrete choice model was created which examined which charging station was selected by EV owners when recharging their vehicle. A MNL model was estimated with the maximum likelihood approach with the BIOGEME package in Python~\citep{Bierlaire2020}, using real charging data for EV owners in the province of Qu\'ebec.

For the error terms for each $i \in N, 1 \leq t \leq T$, the options $j \in C_i^{0t} \cup C_i^{1t}$ are divided into two nests: one for the opt-out option and one for all charging stations. The $\norm{C_i^{0t} \cup C_i^{1t}} \times 2$ factor loading matrix $F$ and $2 \times 2$ diagonal matrix $T$ are given by
\begin{equation}
F = \begin{bmatrix}{cc}
1 & 0 \\
0 & 1 \\
\vdots & \vdots\\
0 & 1
\end{bmatrix}, \quad
T = \begin{bmatrix}{cc}
1 & 0 \\
0 & 1 
\end{bmatrix}.
\end{equation}

\subsection{Distance Dataset}

The user classes, choice sets, and error terms are all identical to the Simple dataset.

The coefficient for distance in the alternative-specific constant has been increased by a factor of ten. More specifically, for each $i \in N, 1 \leq t \leq T$, the alternative-specific constant for each station $j \in C_i^{1t}$ is calculated as
\begin{equation}
\kappa_{ji}^{t} = 1.464 \delta_1 - 0.63 \delta_2 + 0.174 \delta_3,
\end{equation}
with $\delta_1, \delta_2, \delta_3$ defined as in the Simple dataset.

\subsection{Home Charging Dataset}
The set of user classes $N$ includes two user classes for every node in the network: one which has access to home charging, and one which does not. We estimate the  access to home charging via the housing information in the 2016 census~\citep{StatsCan2016}. Based on recommendations from our industrial partners, we assume that 90\% of users in single homes have access to home charging, while 75\% of those in attached homes, and 40\% of those in apartments also have access. The population of each of the two user classes are given by the respective estimates multiplied by a factor of 0.1. 

For user classes $i$ which have access to home charging and for each $1 \leq t \leq T$, the utility for $j \in C_i^{1t}$ is linear in terms of the number of charging outlets, with 
\begin{equation}
\beta_{jik}^{t} = 0.211 k, \forall j \in C_i^{1t}.
\end{equation}
For user classes $i$ which do not have access to home charging and for each $1 \leq t \leq T$, the utility for $j \in C_i^{1t}$ is linear in terms of the number of charging outlets, with 
\begin{equation}
\beta_{jik}^{t} = 0.351 k, \forall j \in C_i^{1t}.
\end{equation}
In both cases, the choice set $C_i^{1t}$ includes all stations which are within ten kilometres of the location of the user class and the alternative-specific constants are identical to the Simple dataset.

For user classes $i$ which do not have access to home charging, the error terms are identical to the Simple dataset. 
For user classes $i$ which have access to home charging and for each $1 \leq t \leq T$, the options $j \in C_i^{0t} \cup C_i^{1t}$ are divided into three nests: one for the opt-out option, one for home charging, and one for all charging stations. The $\norm{C_i^{0t} \cup C_i^{1t}} \times 3$ factor loading matrix $F$ and $3 \times 3$ diagonal matrix $T$ are given by
\begin{equation}
F = \begin{bmatrix}{ccc}
1 & 0 & 0 \\
0 & 1 & 0 \\
0 & 0 & 1 \\
\vdots & \vdots & \vdots\\
0 & 0 & 1 
\end{bmatrix}, \quad
T = \begin{bmatrix}{ccc}
1 & 0 & 0 \\
0 & 1 & 0 \\
0 & 0 & 1
\end{bmatrix}.
\end{equation}

\subsection{LongSpan Dataset}

The mechanisms for the user classes, alternative-specific constants, and error terms are all identical to the Simple dataset. However, the choice sets for each user class now include all stations, not only those within ten kilometres. This, combined with the increased number of stations and the longer time span, results in a significantly more difficult problem to solve.

\subsection{Price Dataset}

The alternative-specific constants, error terms, and choice sets are identical to the LongSpan dataset.

In this dataset, we simulate a price decrease year-by-year, which affects different user classes differently based on their income. The set of user classes $N$ includes five user classes for every node in the network, based on the partitioning in \cite{Javid2016} for income. In the aforementioned work, a logit model for EV acquisition was estimated, with one of the considered factors being the annual household income. The income level was classified as a categorical variable, with the categories defined via income
\begin{itemize}
\item Less than 25 000\$,
\item 25 000\$ - 49 999\$,
\item 50 000\$ - 74 999\$,
\item 75 000 - 99 999\$,
\item Greater or equal to 100 000\$.
\end{itemize} 
In the final estimation of the logit model, the income variable was found to be significantly significant. The utility coefficient for the categorical variable was estimated as 0.443. 

In our work, we estimate the population in each node that falls within each of the five income brackets using the household income field in the 2016 Statistics Canada census~\citep{StatsCan2016}, and assigned each to a user class.\footnote{The census provides data in brackets of 10 000\$, and so the population in certain fields was divided evenly into two user classes (e.g. half of the population of the ``20 000\$ to 29 999\$'' field in the census was assigned to the ``Less than 25 000\$'' user class whereas the other half was added to the ``25 000\$ - 49 999\$'' user class.)} The population of each of the five user classes are given by the respective estimates multiplied by a factor of 0.1, and any user class which would have a population $<1$ are removed. 

An additional term is added to the alternative specific constants for all charging stations based on the income bracket, in increments of 0.443. We then modify the value of the penalisation term each year to account for a decrease in price affecting each user class differently (with the modification affecting the lower income brackets more). 
More specifically, for each $i \in N, 1 \leq t \leq T$, the alternative-specific constant for each station $j \in C_i^{1t}$ is calculated as
\begin{equation}
\kappa_{ji}^{t} = 1.464 \delta_1 - 0.063 \delta_2 + 0.174 \delta_3 + 0.443 \delta_{4i} + 0.443 \left( t - 1 \right) \left( \frac{2 - \delta_{4i}}{4}\right),
\end{equation}
with $\delta_1, \delta_2, \delta_3$ defined as in the Simple dataset and $\delta_{4i}$ given in Table~\ref{TableDelta4i}.

\begin{table}
\centering
\begin{tabular}{lr}
\hline
Income level of user class $i$ & $\delta_{4i}$ \\
\hline
Less than 25 000\$ & -2 \\
25 000\$ - 49 999\$ & -1 \\
50 000\$ - 74 999\$ & 0 \\
75 000 - 99 999\$ & 1 \\
Greater or equal to 100 000\$ & 2
\end{tabular} 
\caption{Values of parameter $\delta_{4i}$}\label{TableDelta4i}
\end{table}

\section{Growth Function model}
\label{AppendixGrowthFunction}
\subsection{Intracity model}
For comparing the GF model of \cite{Anjos2020} to the MC model~\eqref{ModelMaximumCover}, it must be reduced to an intracity form.  More precise definitions and development of each of these variables and equations, we refer to the previous work. Note that some variable names have been changed from the original work to avoid confusion with notation in MC model~\eqref{ModelMaximumCover} and the SL model~\eqref{ModelSinglelevel}. We assume that the city occupies a single urban centre $u$. We also eliminate the path-based constraints from the optimisation model, as these represent users travelling between urban centres. Given these simplifications, we use the following notation:
\begin{itemize}
\item $T$: Set of investment periods.
\item $N$: Set of population centers.
\item $M$: Set of candidate locations.
\item $N_j, j \in M$: Set of locations which are willing to charge at location $j$.
\item $e_j, j \in M$: Maximum number of charging outlets at location $j$. 
\item $r$: Population of the city.
\item $r_i$: Population in location $i$.
\item $l_j, j \in M$: Number of charging outlets already installed at location $j$.
\item $c^U$: Cost for installing a charging outlet at any location.
\item $c^F_j, j \in M$: One-time cost for opening location $j$.
\item $B^t, t \in T$: Budget for year $t$.
\item $\alpha$: Fraction of EV users that choose to charge at home.
\item $a^t, t \in T$: Capacity increase for each charging outlet in year $t$.
\item $S$: Set of segments in the (piecewise linear) GF.
\item $q^{s-1}, q^s, s \in S$: Breakpoints of segment $s$ in the growth function.
\item $m^s, s \in S$: Slope of segment $s$ in the growth function.
\item $o^s, s \in S$: Intercept of segment $s$ in the growth function.
\item $x_j^t, j \in M, t \in T$: Number of charging outlets at station $i$ in year $t$.
\item $y_j^t, j \in M, t \in T$: 1 if station is open in year $t$, 0 otherwise.
\item $w^{st}, s \in S, t \in T$: 1 if the city is at penetration level $s$ at the beginning of year $t$, 0 otherwise.
\item $h_{ij}^t, i \in N,j \in M, t \in T$: Number of EVs based in location $i$ choosing to charge in location $j$ in year $t$.
\item $z^{st}, s \in S, t \in T$: Number of EVs in the city which is at penetration level $s$ at the beginning of year $t$.
\end{itemize}
The list of parameter values can be found in Table~\ref{TableParametersGF}. 

The model used for the comparisons is the following:

\begin{align}
\operatorname{Maximise } & \sum_{j \in M} \sum_{i \in N_j} h_{ij}^{t-1}, \label{MartimModelObjective}\\
\text{subject to}& \sum_{j \in M} c_U \left(x_{j}^ t-x_{j}^{t-1} \right)+\sum_{j \in M} c_j^F \left(y_j^t-y_j^{t-1} \right)  \leq B^t, \quad t \in T, \label{MartimModelBudget}\\
& x_{j}^t \leq e_j y_j^t, \quad j \in M, 1 \leq t \leq T, \label{MartimeModelUpperBound}\\
& x_{j}^t \geq x_{j}^{t-1}, \quad j \in M, 1 \leq t \leq T, \label{MartimModelCantRemoveOutlets}\\
& y_j^t \geq y_j^{t-1}, \, j \in M, 1 \leq t \leq T, \label{MartimModelCantCloseStations}\\
& \sum_{s \in S} z^{st} = \sum_{j \in M} \sum_{i \in N_j} h_{ij}^{t-1}, t \in T, \label{MartimModelSetZ}\\
& q^{s-1}w^{st} \leq z^{st} \leq q^s w^{st}, s \in S, t \in T, \label{MartimModelFindSegment}\\
& \sum_{s \in S} w^{st} \leq 1, t \in T, \label{MartimModelMaxOneSegment} \\
& \sum_{i \in N_j} h_{ij}^t \leq \sum_{i \in N_j} h_{ij}^{t-1} + \frac{r_i}{r} \sum_{s \in S} \left( o^s w^{st}+ \left( m^s-1\right) z^{st} \right), j \in M, t \in T, \label{MartimGrowthFunction}\\
& \sum_{i \in N_j} h_{ij}^{t-1} \leq \sum_{i \in N_j} h_{ij}^{t}, j \in M, t \in T, \label{MartimModelDontLoseEVs}\\
& \alpha \sum_{i \in N_j} h_{ij}^{t} \leq a^t \left( x_j^0+\sum_{t' \leq t} x_j^{t'} \right) j \in N, t \in T. \label{MartimModelCapacity}
\end{align}

The objective function~\eqref{MartimModelObjective} aims to maximise the total number of EV users in the final year.  Constraints~\eqref{MartimModelBudget} are budget constraints, ensuring that the cost of opening charging stations and installing charging outlets does not exceed the budget for that year. Constraints~\eqref{MartimeModelUpperBound} both enforce a maximum number of charging outlets at each station and also ensures that the one-time cost to open charging stations is paid. Constraints~\eqref{MartimModelCantRemoveOutlets} prevent removing charging outlets and  constraints~\eqref{MartimModelCantCloseStations} prevent closing charging stations from one year to the next. Constraints~\eqref{MartimModelSetZ} set the number of EVs at the start of one year as the number at the end of the previous year. Constraints~\eqref{MartimModelFindSegment} find the segment of the growth function that the current EV population is in. Constraints~\eqref{MartimModelMaxOneSegment} ensure that only one segment of the growth function is selected. Constraints~\eqref{MartimGrowthFunction} cap the the number of EVs by the end of the year by following the growth function.   Constraints~\eqref{MartimModelDontLoseEVs} ensure that the total number of EVs does not decrease from year to year. Constraints~\eqref{MartimModelCapacity} are capacity constraints, ensuring that potential new EV users will only decide to purchase an EV if there exists sufficient charging infrastructure.

\begin{table}
\centering
\begin{tabular}{l r}
\hline
Parameter & Value \\
\hline
$T$ & 4 \\
$\abs{N}$ &  317 \\
$\abs{M}$ & 10 \\
$\abs{N_j}$ &  Varies \\
$e_j$ (all stations) &  6 \\ 
$r$ &  181624 \\
$l_j$ (all stations)&  0 \\
$c^U$ (all stations)&  50 \\
$c^F_j$ (all stations)& 100 \\
$B^t$ (per year)& 400 \\
$\alpha$ & 0.566\\
$a^t$ (all years)& $+\infty$ \\
$\abs{S}$ & 5  \\
$q^{s-1}, q^s$ & Varies\\ 
$m^s$ & Varies\\
$o^s$ & Varies \\
\hline
\end{tabular}
\caption{Growth Function parameter values} \label{TableParametersGF}
\end{table}

\subsection{Generating the Growth Function}
The growth function in the GF model gives the number of EVs in the current year as a function of the number of EVs in the previous year. In the absence of the capacity contraints~\eqref{MartimModelCapacity}, the growth function would directly dictate the number of EVs each year via Constraints~\eqref{MartimModelSetZ}. We can ensure that the EV growth remains comparable between the MC and GF models by using the output from the MC model to create the growth function. 

More specifically, we assume there are no EV owners at the start of the optimisation period. While this is not a realistic assumption, it ensures feasibility in the GF model. Given a candidate solution $(\x, \y)$, we solve the MC model~\eqref{ModelMaximumCover} with the desired user classes and parameters over the 20 instances in the dataset. In each instance and for $1 \leq t \leq T$ we calculate the number of users who are covered by $\x$ (given by $\sum_{i \in N} \sum_{r=1}^{R_i} \frac{N_i^{t}}{R_i} w_{i}^{rt}$). To calculate the total number of EVs, we add the new EVs in each year to the EVs from the previous year (or the starting EVs in the case of the first year). We take the average result over all instances for each year as our desired growth function, which mimics perfectly the EV growth from the MC model~\eqref{ModelMaximumCover}. 

\begin{figure}
{\includegraphics[scale=1]{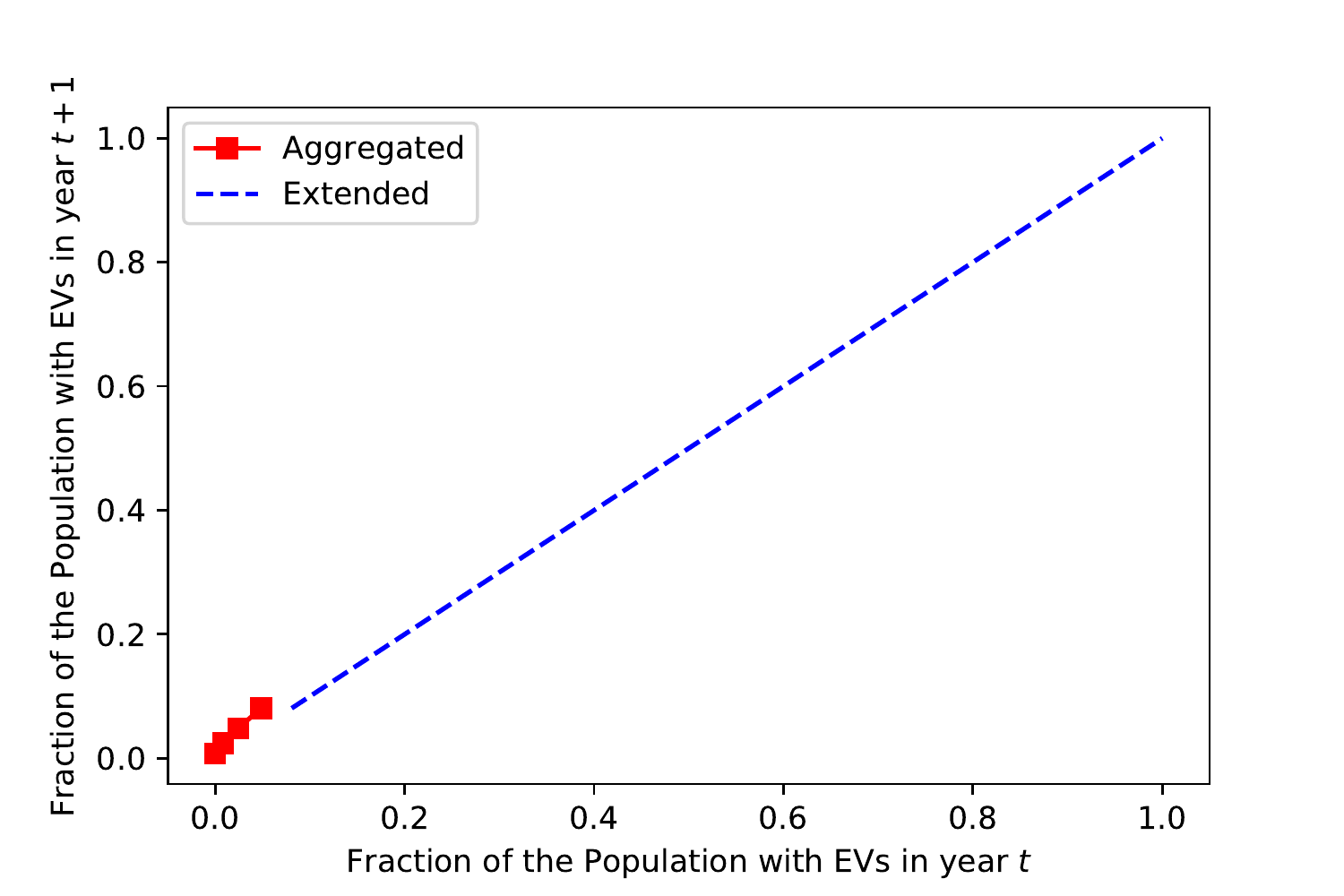}}
\caption{Growth Function Example}
\label{FigureGrowthFunction}
\end{figure}

After normalising for the population--which gives the percentage of the population with EVs in the following year given the percentage of the population with EVs in the current year--, we extend the growth function to cover the entire $[0,1]$ domain. Both of these steps allow the growth function to be used regardless of population. An example of the normalised, extended growth function is given in Figure~\ref{FigureGrowthFunction}. 

\end{document}